\documentclass[11pt]{article}
\usepackage[T1]{fontenc}
\usepackage[utf8]{inputenc}
\usepackage{amsmath,amssymb,amsthm}
\usepackage{listings}
\usepackage[margin=1.15in]{geometry}
\usepackage{hyperref}

\lstdefinestyle{python}{
  language=Python,
  basicstyle=\ttfamily\scriptsize,
  numbers=none,
  breaklines=true,
  showstringspaces=false,
  keepspaces=true,
  columns=fullflexible,
  frame=single,
  framesep=3pt,
  upquote=true,
  literate={-}{-}1
}

\lstdefinestyle{grid}{
  basicstyle=\ttfamily\scriptsize,
  numbers=none,
  frame=none,
  columns=fullflexible,
  keepspaces=true,
  breaklines=true,
  showstringspaces=false
}

\theoremstyle{plain}
\newtheorem{theorem}{Theorem}
\newtheorem{lemma}[theorem]{Lemma}
\newtheorem{proposition}[theorem]{Proposition}
\newtheorem{corollary}[theorem]{Corollary}
\theoremstyle{definition}
\newtheorem{definition}[theorem]{Definition}
\newtheorem{remark}[theorem]{Remark}

\DeclareMathOperator{\slack}{sl}
\DeclareMathOperator{\comp}{comp}
\newcommand{\Z}{\mathbb{Z}}

\newcommand{\Rmax}{\overline{\Z}}

\title{The maximum area of a snake polyomino in a $12 \times n$ rectangle}
\author{Tong Niu}
\date{\today}

\begin{document}
\maketitle

\begin{abstract}
A \emph{snake polyomino} is a set of cells of the square lattice whose
cells, under edge adjacency, induce a path. Let $T(m,n)$ be the largest
number of cells of a snake polyomino contained in an $m \times n$
rectangle, equivalently the largest number of vertices of an induced path
in the grid graph $P_m \times P_n$ (OEIS \texttt{A360917}). Blondin
Mass\'e and Goupil determined $T(m,n)$ for all $m \le 5$ and proposed the
general case as a conjecture; for $m = 12$, entry \texttt{A397153} of the
OEIS records a conjectured generating function, contributed by Goupil in
July 2026, which has remained unproved. We prove it. For every $n \ge 1$
other than $n \in \{2,3,\dots,10,12\}$ one has $T(12,n) = 44k + c_r$ where
$n = 6k+r$ with $0 \le r < 6$ and $(c_0,\dots,c_5) = (5,12,19,27,34,42)$,
and the ten remaining values are $18,26,33,40,47,55,62,70,77$ and $92$.
The proof has two independent halves. The upper bound comes from a
column-profile transfer system read in the max-plus semiring: a single
entrywise inequality $v_{46} \le 44 + v_{40}$ between two value vectors on
its useful core propagates to all larger indices because the max-plus
action is monotone and additively homogeneous. The lower bound comes from
explicit periodic snakes, one for each residue class, which are checked
against the definition of a snake and not against the transfer system. The
same method settles every height up to $12$; heights $6$ through $11$ had
no formula in the literature, and we record them as well.
\end{abstract}



\section{Introduction}

\subsection{The problem}

Write $[m] = \{1, \dots, m\}$ and let $G_{m,n}$ be the grid graph on the
cell set $[m] \times [n]$, two cells being adjacent when they agree in one
coordinate and differ by one in the other. A nonempty $S \subseteq [m]
\times [n]$ is a \emph{snake polyomino} when the subgraph $G_{m,n}[S]$
induced by $S$ is a path; the single cell counts as a path on one vertex.
Equivalently $S$ is connected, contains no cycle, and no cell of $S$ has
three neighbours in $S$. Put
\[
  T(m,n) \;=\; \max \{\, |S| \;:\; S \subseteq [m] \times [n]
  \text{ a snake polyomino} \,\}.
\]
Snake polyominoes were introduced in this form by Blondin Mass\'e and
Goupil \cite{BG2026}, who call them the two-dimensional binary words whose
cell graph is a chain; the quantity $T(m,n)$ is the maximum number of
vertices of an induced path, sometimes called a \emph{snake} or an induced
path of maximum order, in $P_m \times P_n$. The table of values is entry
\texttt{A360917} of the On-Line Encyclopedia of Integer Sequences
\cite{A360917}.

Deciding the maximum order of an induced path is NP-hard for general
graphs \cite{GareyJohnson}, so no generic tool settles $T(m,n)$. Beluhov
\cite{Beluhov2023} did prove the asymptotic
$T(m,n) = \tfrac{2}{3}mn + O(m+n)$, but an asymptotic is not a table of
values, and the exact values have been reached one height at a time.
Blondin Mass\'e and Goupil \cite{BG2026} give exact formulas for
$m \in \{1,2,3,4,5\}$ and write that ``our underlying conjecture is that
there exists an explicit expression for the maximal area of snake-like
polyominoes in any given rectangle''. Their companion paper with Gauthier
and Toure \cite{GGT2026} puts it plainly: ``the problem of finding the
length of maximal snake polyominoes in an arbitrary rectangle is still
open''. The best exact method available is complete enumeration, which
they push to $15 \times 15$.

Our subject is the column $m = 12$. Goupil contributed it to the OEIS in
July 2026 as entry \texttt{A397153} \cite{A397153}, whose only formula line
reads, verbatim,
\begin{quote}
Conjectured g.f.: \texttt{x*(12 + 6*x + 8*x\^{}2 + 7*x\^{}3 + 7*x\^{}4 +
7*x\^{}5 - 4*x\^{}6 + x\^{}7 + 2*x\^{}10 - x\^{}11 - x\^{}16 + x\^{}17 -
x\^{}18)/((1 - x)\^{}2*(1 + x)*(1 - x + x\^{}2)*(1 + x + x\^{}2))}.
\end{quote}
Howroyd extended the data to $a(60)$ by brute force in July 2026. An
editorial pass in August 2026 then removed a second conjectured closed
form while deliberately keeping this one, on the ground that it shows the
structure better. The word \emph{Conjectured} has stood ever since.

\subsection{What is proved}

\begin{theorem}\label{thm:main}
Put $(c_0,c_1,c_2,c_3,c_4,c_5) = (5,12,19,27,34,42)$ and, for $n \ge 1$,
write $n = 6k+r$ with $k \ge 0$ and $0 \le r < 6$. Then
\[
  T(12,n) \;=\; 44k + c_r
  \qquad\text{for } n = 1,\ n = 11, \text{ and every } n \ge 13,
\]
while the ten remaining values are
\[
  \bigl(T(12,2), \dots, T(12,10)\bigr) = (18,26,33,40,47,55,62,70,77),
  \qquad T(12,12) = 92 .
\]
Equivalently
\[
  \sum_{n \ge 1} T(12,n)\, x^n
  \;=\; \frac{x\,N(x)}{(1-x)^2 (1+x)(1-x+x^2)(1+x+x^2)},
\]
where
\begin{align*}
  N(x) \;=\; 12 + 6x + 8x^2 + 7x^3 &+ 7x^4 + 7x^5 - 4x^6 + x^7 \\
       &+ 2x^{10} - x^{11} - x^{16} + x^{17} - x^{18},
\end{align*}
which is Goupil's conjectured generating function for \texttt{A397153}.
\end{theorem}

The denominator factors as $(1-x)^2(1+x)(1-x+x^2)(1+x+x^2) = (1-x)(1-x^6)$,
so the shape of the answer --- linear growth with a period-$6$ correction
--- is already visible in the conjecture. What the conjecture does not say
is that the period really is $6$ and the shift really is $44$. That is the
part that has to be proved, and sixty terms cannot do it: they cannot tell
a genuine period from an artefact of the range they happen to cover.

In particular the growth rate is $22/3$ cells per column. That number is
what makes the problem resist the machinery already in print. A snake
polyomino is, after all, a set of cells of maximum degree at most two, and
Blondin Mass\'e, Goupil, L'Heureux and Marin \cite{BGLM2025} compute the
exact maximum area of such a set in an $h \times w$ rectangle; at $h = 12$
their bound grows at $8$ cells per column, strictly larger than $22/3$.
Goupil observed on \texttt{A360917} that at height four the degree-two
bound is attained, and that is how the small heights were reachable. At
height twelve it is not attained, and the miss is $2/3$ of a cell per
column, so a genuinely different upper-bound argument is needed. In the
other direction the obvious construction --- alternate full columns of
twelve cells with a single connector cell at each turn, which is an
induced path --- gives only $13/2$ cells per column. The truth therefore
sits strictly between the two elementary bounds,
\[
  \tfrac{13}{2} \;=\; 6.5 \;<\; \tfrac{22}{3} \;=\; 7.33\ldots \;<\; 8 ,
\]
and neither of them can decide it.

The method also settles every smaller height. No height above $5$ has a
formula in the literature --- \cite{BG2026} stops at $m = 5$ --- and for
heights $7$ to $11$ there is not even an entry in the OEIS: a search on
the first twelve terms of each of those five sequences returns nothing. So
the closed forms recorded here appear to be the first. Height $6$ is entry
\texttt{A397757} \cite{A397757}, where a comment of Howroyd from July 2026
states that the formula there ``is correct. It was confirmed via the
transfer-matrix method''. No proof was published; our certificate turns
the claim into a verifiable statement.

\begin{theorem}\label{thm:all}
For each $m \le 12$ there are a period $P_m$, a shift $C_m$, an onset
$n_m$ and constants $c^{(m)}_0, \dots, c^{(m)}_{P_m - 1}$ --- the first
three in Table~\ref{tab:all}, the constants in Table~\ref{tab:const} ---
such that
\[
  T(m,n) \;=\; C_m k + c^{(m)}_r
  \qquad\text{whenever } n = P_m k + r,\ 0 \le r < P_m,\ n \ge n_m,
\]
and the finitely many values $T(m,n)$ with $n < n_m$ are those listed in
Table~\ref{tab:small}. The growth rate is $C_m / P_m$ cells per column.
\end{theorem}

The periods in Table~\ref{tab:all} are not a pattern one would guess:
they are $1,2,1,3,7,4,6,8,4,18,8,6$ for $m = 1, \dots, 12$, and the
quasi-linear regime at height $11$ only begins at $n = 52$. That is why
we insist on a certificate instead of an extrapolation from the first
twenty terms.

\subsection{Method, and what is load-bearing}

Section~\ref{sec:machine} sets up a transfer system that reads a snake one
column at a time. Roughly, its states record, for the frontier column,
which cells are occupied, how much room each of them still has, and which
of them lie on a common component of the part already read. The one fact we need about
it is \emph{completeness}: every snake of column span $1..n$ traces a
legal trajectory of length $n$ whose weights sum to the number of its
cells and whose last configuration has a single component
(Lemma~\ref{lem:complete}). The converse, that every trajectory comes from
a snake, is Lemma~\ref{lem:sound}; together they make the computed values
exactly $T(m,n)$. Only completeness is load-bearing for
Theorem~\ref{thm:main}, since there the lower bound is carried by snakes
written down explicitly. The point is that a fault in the machine can then
only \emph{weaken} the upper bound. For Theorem~\ref{thm:all} we use
both.

Section~\ref{sec:upper} reads the system in the max-plus semiring. Writing
$v_j$ for the vector of best weights of trajectories of length $j$, the
recursion is $v_{j+1} = M \otimes v_j$ with $M$ the transfer matrix. The
max-plus action is monotone and additively homogeneous, so from a single
entrywise inequality
\[
  v_{a + p} \;\le\; c + v_a
\]
one gets $v_{j+p} \le c + v_j$ for every $j \ge a$ by applying $M$
repeatedly, hence $T(m,n+p) \le c + T(m,n)$ for all large $n$. For
$m = 12$ the inequality holds with $p = 6$, $c = 44$ and $a = 40$, over a
useful core of 173\,994 states. This is a weaker and much more
forgiving object than the exact periodicity $v_{a+p} = c + v_a$ of the
value vector that one might reach for first. That identity needs the whole
state space to have entered its periodic regime, and it already fails at
height $6$, where the correction term has period $8$ although
$T(6,\cdot)$ has period $4$. The inequality, on the other hand, is
satisfied for free by every state whose value grows slower than $c/p$.
Certificates of this shape for fixed-width transfer computations are
standard --- they are how the domination numbers of grid
graphs were settled \cite{GPRT2011,Alanko2011} --- and the contribution
here is the application, not the technique.

Section~\ref{sec:lower} gives the matching lower bound by explicit
periodic snakes $S_k = A \cdot B^k \cdot C$, one family per residue class.
Their verification uses only the definition of a snake and three finite
checks, and in particular it does not use the transfer system. The
criterion is: maximum degree at most two, connected, and some cell of
degree at most one. We insist on this exact form because the shorter
criterion ``maximum degree at most two and $|E| = |S| - 1$'' is false. For a
graph of maximum degree at most two, $|E| = |S| - (\text{number of
\emph{path} components})$, so cycle components contribute no edge deficit
at all; a domino placed far from an induced $4$-cycle has six cells, five
edges and maximum degree two, and is not a snake.

Section~\ref{sec:verify} reports the verification. The computation was
carried out twice, by two implementations written independently of each
other; they agree on the state count, the edge count, the number of
accepting states and the multiset of edge weights at every height from $1$
to $9$, and on the useful core itself at heights $10$, $11$ and $12$. The
completeness lemma was tested against \emph{every} snake polyomino of
fifteen small rectangles. And the values reproduce all $60$ published
terms of \texttt{A397153}, all $10000$ terms of the b-file of
\texttt{A397757}, the whole available part of the b-file of
\texttt{A360917}, together with a brute-force search over induced paths.

\section{The transfer system}\label{sec:machine}

Throughout this section $m \ge 1$ is fixed. For a snake polyomino
$S \subseteq [m] \times [n]$ and $1 \le j \le n$ write
$S_j = S \cap ([m] \times [j])$ for its \emph{prefix} and $G_j$ for the
subgraph of $G_{m,n}$ induced by $S_j$. Cells are pairs $(i,j)$ with $i$
the row and $j$ the column.

\subsection{Prefixes of a snake}

\begin{lemma}\label{lem:span}
Let $S \subseteq [m] \times [n]$ be a snake polyomino. Then the set of
columns it meets is an interval of integers. Consequently, putting
\[
  f(j) \;=\; \max\{\, |S| : S \subseteq [m]\times[j] \text{ a snake
  meeting column } 1 \text{ and column } j \,\},
\]
we have $T(m,n) = \max_{1 \le j \le n} f(j)$.
\end{lemma}

\begin{proof}
If column $j$ is empty while columns $j' < j$ and $j'' > j$ are not, then
every path of $G_{m,n}$ from a cell of column $j'$ to a cell of column
$j''$ must use a cell of column $j$, so $G_{m,n}[S]$ is disconnected,
which a path is not. For the second claim, a snake $S$ in $[m]\times[n]$
meets exactly the columns $u, u+1, \dots, w$; translating it $u-1$ steps
to the left gives a snake of the same size in $[m] \times [w-u+1]$ meeting
the first and the last column, so $T(m,n) \le \max_{j \le n} f(j)$.
Conversely a snake counted by $f(j)$ for $j \le n$ is a snake in
$[m]\times[n]$.
\end{proof}

For $v \in S_j$ put $\slack_j(v) = 2 - \deg_{G_j}(v)$; this counts the
neighbours $v$ may still acquire without breaking the degree bound.

\begin{lemma}\label{lem:prefix}
$G_j$ is a disjoint union of paths, and every component $C$ of $G_j$
satisfies $\sum_{v \in C} \slack_j(v) = 2$.
\end{lemma}

\begin{proof}
$G_j$ is an induced subgraph of the path $G_{m,n}[S]$, and an induced
subgraph of a path is a disjoint union of paths. A path on $t$ vertices
has $t-1$ edges, so its degrees sum to $2t-2$ and its slacks to
$2t - (2t-2) = 2$.
\end{proof}

Now fix $S$ and $j$ with $S_j \ne \emptyset$ and set
\begin{align*}
  O_j &= \{\, i : (i,j) \in S \,\}, &
  s_j(i) &= \slack_j\bigl((i,j)\bigr) \quad (i \in O_j), \\
  O_j^{+} &= \{\, i \in O_j : s_j(i) \ge 1 \,\}, &
  c_j &= \text{number of components of } G_j,
\end{align*}
let $P_j$ be the partition of $O_j^{+}$ in which $i$ and $i'$ lie in one
block exactly when $(i,j)$ and $(i',j)$ lie in one component of $G_j$, and
put
\[
  \delta_j \;=\; 2 c_j - \sum_{i \in O_j} s_j(i).
\]
Cells of $O_j^{+}$ are called \emph{open}. By Lemma~\ref{lem:prefix} an
open cell has degree at most $1$ in $G_j$, so it is an endpoint of its
component; a component of $G_j$ therefore owns at most two open cells of
column $j$, and the slacks of the open cells of one component sum to at
most $2$.

\begin{lemma}[dead slack]\label{lem:dead}
$\delta_j = \sum \{\, \slack_j(v) : v \in S_j,\ v \notin [m]\times\{j\}
\,\}$, and $0 \le \delta_1 \le \delta_2 \le \cdots \le \delta_n \le 2$.
\end{lemma}

\begin{proof}
The displayed identity is Lemma~\ref{lem:prefix} summed over the $c_j$
components of $G_j$, minus the slack carried by column $j$. A cell in a
column smaller than $j$ has no neighbour in column $j+1$, so its degree,
and therefore its slack, is the same in $G_j$ and in $G_{j+1}$; hence
\[
  \delta_{j+1} \;=\; \delta_j + \sum_{i \in O_j}
  \slack_{j+1}\bigl((i,j)\bigr) \;\ge\; \delta_j .
\]
Finally $S_n = S$ is a single path, so $\delta_n = 2 - \sum_{i \in O_n}
s_n(i) \le 2$ by Lemma~\ref{lem:prefix}, and $\delta_1 \ge 0$ because
slacks are nonnegative.
\end{proof}

\begin{lemma}\label{lem:open}
Let $1 \le j < n$ and let $C$ be a component of $G_j$. Then
\begin{enumerate}
\item[(a)] $C$ owns a cell $(i,j)$ with $(i,j+1) \in S$, and every such
  cell is open;
\item[(b)] in particular $C$ owns an open cell of column $j$, so
  $c_j = |P_j|$ and $\delta_j = 2|P_j| - \sum_{i \in O_j} s_j(i)$.
\end{enumerate}
\end{lemma}

\begin{proof}
Suppose no cell of $C$ has a neighbour outside $S_j$. Every neighbour in
$S$ of a cell of $C$ lies in $S_j$, hence in $C$, so $C$ is a connected
component of $G_{m,n}[S]$; as that graph is connected, $C = S$, and then
$S$ misses column $n > j$, contradicting Lemma~\ref{lem:span}. So some
cell of $C$ has a neighbour $v \in S \setminus S_j$; grid edges change one
coordinate by one, so that cell is some $(i,j) \in C$ and $v = (i,j+1)$.
Such a cell has a neighbour outside $S_j$ and degree at most $2$ in
$G_{m,n}[S]$, hence degree at most $1$ in $G_j$, i.e.\ $s_j(i) \ge 1$.
Part (b) follows because distinct components give distinct blocks and
every component gives a nonempty one.
\end{proof}

\subsection{Profiles and transitions}

\begin{definition}
A \emph{profile} of height $m$ is a triple $\pi = (O, s, P)$ where
$\emptyset \ne O \subseteq [m]$, $s : O \to \{0,1,2\}$, and $P$ is a
partition of $O^{+} := \{i \in O : s(i) \ge 1\}$. We write
$\delta(\pi) = 2|P| - \sum_{i \in O} s(i)$.
\end{definition}

So Lemma~\ref{lem:open}(b) says that for $j < n$ the triple
\[
  \sigma(S_j) \;:=\; \pi_j \;:=\; (O_j, s_j, P_j)
\]
is a profile with $\delta(\pi_j) = \delta_j$; we call it the \emph{profile
of the prefix} $S_j$.

\begin{definition}[initial profiles]
For $\emptyset \ne O \subseteq [m]$ let $\iota(O) = (O, s, P)$ where
$s(i) = 2 - |\{i-1,i+1\} \cap O|$ and $P$ is the partition of $O^{+}$ in
which $i < i'$ share a block exactly when $\{i, i+1, \dots, i'\}
\subseteq O$.
\end{definition}

A direct count gives $\delta(\iota(O)) = 0$. Indeed, a maximal run of $O$
of length $1$ carries slack $2$ on its single cell, while a run of length
$t \ge 2$ carries slack $1$ on each of its two extreme cells and $0$
inside, so each of the $|P|$ runs contributes exactly $2$.

\begin{definition}[transition]\label{def:trans}
Let $\pi = (O,s,P)$ be a profile and $\emptyset \ne O' \subseteq [m]$. Let
$H(\pi, O')$ be the graph with vertex set $P \sqcup O'$ (the blocks of
$\pi$ and the rows of $O'$) and edges
\[
  \{i, i+1\} \ \text{ for } i, i+1 \in O',
  \qquad
  \{i, B\} \ \text{ for } i \in O \cap O' \text{ and } i \in B \in P,
\]
and for $i \in O'$ put
\[
  d'(i) = [\,i \in O\,] + |\{i-1,i+1\} \cap O'| , \qquad s'(i) = 2 - d'(i).
\]
The pair $(\pi, O')$ is \emph{admissible} when
\begin{enumerate}
\item[(T1)] $s(i) \ge 1$ for every $i \in O \cap O'$;
\item[(T2)] $d'(i) \le 2$ for every $i \in O'$;
\item[(T3)] $H(\pi, O')$ is a forest;
\item[(T4)] every block $B \in P$ meets an edge of $H(\pi, O')$.
\end{enumerate}
For an admissible pair let $\comp'(\pi,O')$ be the number of connected
components of $H(\pi,O')$ and
$\delta'(\pi,O') = 2\comp'(\pi,O') - \sum_{i \in O'} s'(i)$. The pair is
\emph{live} when $\delta'(\pi,O') \le 2$ and every component of
$H(\pi,O')$ contains a row $i \in O'$ with $s'(i) \ge 1$; its
\emph{successor} is then the profile
$\tau(\pi, O') = \bigl(O', s', P'\bigr)$, where $P'$ is the partition of
$\{i \in O' : s'(i) \ge 1\}$ by the components of $H(\pi, O')$.
\end{definition}

Condition (T2) forces $s'(i) \ge 0$, and it is purely local: $d'(i) = 3$
requires $i \in O$ together with $i-1 \in O'$ and $i+1 \in O'$, so (T2)
fails exactly when some $i$ has $i-1, i, i+1 \in O'$ and $i \in O$. One
test on three consecutive rows, then. The only edges of $H(\pi,O')$
incident with a block are the horizontal ones, so (T4) says simply that
every block of $P$ contains a row of $O \cap O'$. Three counting
identities make the rest cheap; they are also what the two implementations
use as filters.

\begin{lemma}\label{lem:counting}
Let $(\pi, O')$ satisfy (T3), and put $r = |O \cap O'|$ and
$u = |\{ i : i, i+1 \in O' \}|$. Then
\[
  \comp'(\pi, O') = |P| + |O'| - u - r ,
  \qquad
  \sum_{i \in O'} s'(i) = 2|O'| - r - 2u ,
  \qquad
  \delta'(\pi, O') = 2|P| - r .
\]
\end{lemma}

\begin{proof}
$H(\pi,O')$ has $|P| + |O'|$ vertices and $u + r$ edges, and a forest has
(vertices $-$ edges) components, which is the first identity. Summing
$d'$ over $O'$ counts each horizontal edge once and each vertical edge
twice, giving $\sum_{i} d'(i) = r + 2u$ and hence the second. Substituting
both into the definition of $\delta'$ gives
$2(|P| + |O'| - u - r) - (2|O'| - r - 2u) = 2|P| - r$.
\end{proof}

So $\comp' = 1$ reads $|P| + |O'| = u + r + 1$, and the inequality
$\delta' \le 2$ reads $r \ge 2|P| - 2$. A third reading of the same
identities, $u + r \le |P| + |O'| - 1$, is necessary for (T3).

\begin{definition}[the machine]\label{def:machine}
Let $Q = Q_m$ be the smallest set of profiles of height $m$ that contains
$\iota(O)$ for every $\emptyset \ne O \subseteq [m]$ and is closed under
$\tau$ on live admissible pairs. Adjoin one further symbol $\top$ and put
$\widehat{Q} = Q \cup \{\top\}$. The \emph{edges} are
\begin{itemize}
\item $\pi \to \tau(\pi, O')$ of weight $|O'|$, for $\pi \in Q$ and
  $(\pi,O')$ live admissible;
\item $\pi \to \top$ of weight
  $w_\top(\pi) = \max \{\, |O'| : (\pi,O') \text{ admissible with }
  \comp'(\pi,O') = 1 \,\}$, when that set is nonempty;
\end{itemize}
and $\top$ has no outgoing edge. The \emph{accepting} states are
$\mathrm{Acc} = \{\, \pi \in Q : |P(\pi)| = 1 \,\} \cup \{\top\}$.
\end{definition}

The symbol $\top$ is not a convenience. A snake can end with its last
column \emph{saturated}: every cell of the last column already of degree
$2$, so the single component owns no open cell there and
Lemma~\ref{lem:open}(b) --- which assumed $j < n$ --- does not apply. Such
a configuration is a perfectly good snake that cannot be continued, and so
is not a live state; $\top$ is where those terminal configurations get
collected. We remark that leaving them out is not a harmless
simplification: it makes the computed value of $T(3,6)$ come out as $13$
instead of the true $14$.

\subsection{Completeness}

\begin{lemma}[completeness]\label{lem:complete}
Let $S \subseteq [m]\times[n]$ be a snake polyomino meeting column $1$ and
column $n$, with profiles $\pi_j = (O_j, s_j, P_j)$ for $1 \le j \le n$.
Then $\pi_1 = \iota(O_1)$; for every $j < n$ the pair $(\pi_j, O_{j+1})$
is admissible with $\comp'(\pi_j, O_{j+1}) = c_{j+1}$; for every
$j < n-1$ it is live with $\tau(\pi_j, O_{j+1}) = \pi_{j+1}$; and
$c_n = 1$. Consequently $\pi_1, \dots, \pi_{n-1} \in Q$ and there is a
walk
\[
  \pi_1 \longrightarrow \pi_2 \longrightarrow \cdots
  \longrightarrow \pi_{n-1} \longrightarrow x ,
  \qquad x \in \mathrm{Acc},
\]
(of length $n-1$; for $n = 1$ the walk is the single accepting state
$\pi_1$) along which the weights, together with $|O_1|$, sum to at least
$|S|$.
\end{lemma}

\begin{proof}
In column $1$ a cell has neighbours in $S_1$ only above and below it, so
$s_1(i) = 2 - |\{i-1,i+1\} \cap O_1|$, and the components of $G_1$ are the
maximal runs of $O_1$; that is $\pi_1 = \iota(O_1)$.

Fix $j < n$ and write $H = H(\pi_j, O_{j+1})$. (T1) holds because
$i \in O_j \cap O_{j+1}$ means the cell $(i,j)$ has the neighbour
$(i,j+1)$ in $S$, so by Lemma~\ref{lem:open}(a) it is open. (T2) holds
because $d'(i) = \deg_{G_{j+1}}((i,j+1)) \le 2$, as $S$ has maximum degree
$2$.

For (T3) and the component count, consider the map $\varphi$ on the
vertices of $G_{j+1}$ that sends a cell $(i,j+1)$ of column $j+1$ to the
row $i \in O_{j+1}$, and a cell $v \in S_j$ to the block
$C \cap O_j^{+} \in P_j$, where $C$ is the component of $G_j$ containing
$v$; the block is nonempty by Lemma~\ref{lem:open}(b). Under $\varphi$ the
edges of $G_{j+1}$ inside $S_j$ become loops, the vertical edges of column
$j+1$ become the vertical edges of $H$, and a horizontal edge
$\{(i,j),(i,j+1)\}$ becomes $\{i, B\}$ with $i \in B$; distinct
horizontal edges give distinct edges of $H$ because their rows differ. So
$H$ is obtained from the forest $G_{j+1}$ (Lemma~\ref{lem:prefix}) by
contracting each component of $G_j$ to a single vertex, and no parallel
edges arise. Contracting disjoint connected subgraphs changes neither the
cycle rank nor the number of components: the components of $G_j$ are
trees, so contracting them deletes exactly $|S_j| - c_j$ edges and removes
exactly $|S_j| - c_j$ vertices, leaving $|E| - |V|$ and the component count
untouched. Since $G_{j+1}$ is a forest its cycle rank is $0$, hence so is
that of $H$; that is (T3), and $\comp'(\pi_j, O_{j+1}) = c_{j+1}$.

For (T4), let $B = C \cap O_j^{+}$ with $C$ a component of $G_j$. By
Lemma~\ref{lem:open}(a) some $(i,j) \in C$ has $(i,j+1) \in S$ and is
open, so $i \in B \cap O_j \cap O_{j+1}$ and $\{i, B\}$ is an edge of $H$.

Now let $j < n-1$. Then $\delta'(\pi_j, O_{j+1}) = \delta_{j+1} \le 2$ by
Lemmas~\ref{lem:dead} and~\ref{lem:open}(b), and every component of $H$
corresponds to a component of $G_{j+1}$, which owns an open cell of column
$j+1$ by Lemma~\ref{lem:open}(b) applied at $j+1 < n$. So the pair is
live. Its successor is $(O_{j+1}, s', P')$ with
$s'(i) = 2 - \deg_{G_{j+1}}((i,j+1)) = s_{j+1}(i)$ and with $P'$ the
partition by components of $H$, which is $P_{j+1}$; that is
$\tau(\pi_j, O_{j+1}) = \pi_{j+1}$.

Finally $G_n = G_{m,n}[S]$ is a path, so $c_n = 1$. If the pair
$(\pi_{n-1}, O_n)$ is live then $\pi_n = \tau(\pi_{n-1},O_n) \in Q$ and
$|P_n| = c_n = 1$, so $\pi_n \in \mathrm{Acc}$ and we may take
$x = \pi_n$. If it is not live then, $\comp'$ being $1$, the pair is
admissible with $\comp' = 1$, so $w_\top(\pi_{n-1}) \ge |O_n|$ and we may
take $x = \top$. In both cases the weights along the walk are $|O_2|,
\dots, |O_{n-1}|$ followed by a last weight at least $|O_n|$, and
$|S| = \sum_{j=1}^{n} |O_j|$.
\end{proof}

\subsection{Soundness}

The converse of Lemma~\ref{lem:complete} also holds. We prove it because
it is what makes the computed values \emph{equal} to $T(m,n)$ rather than
merely an upper bound, and that equality is convenient for the small
widths at the heights treated in Theorem~\ref{thm:all}. For
Theorem~\ref{thm:main} it is \emph{not} used at all: there the upper bound
comes from Lemma~\ref{lem:complete} alone and the lower bound from snakes
written down explicitly, so a fault in the machine could not produce the
headline result.

Note first that an edge between two profiles determines its own weight:
the edge $\pi \to \tau(\pi,O')$ has weight $|O'|$ and $O'$ is the first
component of $\tau(\pi,O')$.

\begin{lemma}[soundness]\label{lem:sound}
Let $\pi_1 \to \cdots \to \pi_{n-1} \to x$ be a walk in the machine with
$\pi_1 = \iota(O_1)$ an initial profile and $x \in \mathrm{Acc}$. Put
$O_j = O(\pi_j)$ for $j \le n-1$, and let $O_n = O(x)$ if $x \in Q$, while
if $x = \top$ let $O_n$ be any set with $(\pi_{n-1}, O_n)$ admissible,
$\comp'(\pi_{n-1},O_n) = 1$ and $|O_n| = w_\top(\pi_{n-1})$. Then
\[
  S \;=\; \bigcup_{j=1}^{n} O_j \times \{j\}
\]
is a snake polyomino in $[m]\times[n]$ meeting column $1$ and column $n$,
of size $|S| = \sum_{j} |O_j|$, which is the accumulated weight of the
walk; and $\pi_j = \sigma(S_j)$ for every $j \le n-1$.
\end{lemma}

\begin{proof}
We show by induction on $j \le n-1$ that $G_j$ is a forest of maximum
degree at most $2$ with exactly $|P_j|$ components, that every component
owns an open cell of column $j$, and that $\pi_j$ is the profile of $S_j$.
For $j = 1$ this is the computation that gave
$\delta(\iota(O_1)) = 0$: the components of $G_1$ are the maximal runs of
$O_1$, each a path, with the slacks and the partition prescribed by
$\iota$.

Assume it at $j$, and let $O' = O_{j+1}$, $H = H(\pi_j, O')$,
$r = |O_j \cap O'|$ and $u = |\{i : i, i+1 \in O'\}|$. For degrees: a cell
$(i,j+1)$ has $\deg_{G_{j+1}} = d'(i) \le 2$ by (T2); a cell $(i,j)$ gains
a neighbour exactly when $i \in O'$, and then (T1) gives
$s_j(i) \ge 1$, i.e.\ its old degree was at most $1$; cells in columns
below $j$ gain nothing. So $G_{j+1}$ has maximum degree at most $2$.

By the induction hypothesis the blocks of $P_j$ are precisely the open
cells of column $j$ grouped by the components of $G_j$, so the map
$\varphi$ of the proof of Lemma~\ref{lem:complete} is available and
exhibits $H$ as $G_{j+1}$ with every component of $G_j$ contracted to a
point. Contraction of connected subgraphs does not change the number of
connected components, so $G_{j+1}$ and $H$ have the same number, namely
$\comp'(\pi_j,O')$. On the other hand $G_j$ is a forest with $|P_j|$
components, so it has $|S_j| - |P_j|$ edges, and therefore
\[
  |V(G_{j+1})| - |E(G_{j+1})|
  = \bigl(|S_j| + |O'|\bigr) - \bigl(|S_j| - |P_j| + u + r\bigr)
  = |P_j| + |O'| - u - r
  = \comp'(\pi_j, O'),
\]
the last equality being Lemma~\ref{lem:counting}, which applies because
(T3) holds. A graph is a forest exactly when its vertex count minus its
edge count equals its number of components, so $G_{j+1}$ is a forest.
Finally $s'(i) = 2 - \deg_{G_{j+1}}((i,j+1))$ and $P'$ is the partition of
the open cells of column $j+1$ by the components of $H$, hence by the
components of $G_{j+1}$; so $\pi_{j+1} = \tau(\pi_j,O') = \sigma(S_{j+1})$,
and liveness gives that each component owns an open cell. That closes the
induction.

For the last column, suppose first $x = \pi_n \in Q$. The same degree and
forest argument applies to the step $j = n-1$, so $G_n$ is a forest of
maximum degree at most $2$ with $|P_n|$ components, and $|P_n| = 1$
because $x \in \mathrm{Acc}$. A connected forest of maximum degree at
most $2$ is a path. If instead $x = \top$, the chosen $O_n$ makes
$(\pi_{n-1}, O_n)$ admissible with $\comp' = 1$, and the same argument
gives that $G_n$ is a forest of maximum degree at most $2$ with one
component, again a path. In both cases $O_1$ and $O_n$ are nonempty, so
$S$ meets the first and the last column.
\end{proof}

\begin{corollary}\label{cor:exact}
$f(j) = F(j)$ for every $j \ge 1$, where
$F(j) = \max_{x \in \mathrm{Acc}} v_j[x]$ as in Section~\ref{sec:upper};
consequently $T(m,n) = G(n) = \max_{j \le n} F(j)$ for every $n \ge 1$.
\end{corollary}

\begin{proof}
Lemma~\ref{lem:complete} gives $f(j) \le F(j)$
(Proposition~\ref{prop:upper-basic}). Conversely a walk realising
$v_j[x]$ for $x \in \mathrm{Acc}$ yields, by Lemma~\ref{lem:sound}, a
snake of column span $1..j$ whose size is the accumulated weight, so
$f(j) \ge F(j)$. The second statement is Lemma~\ref{lem:span}.
\end{proof}

\section{The upper bound}\label{sec:upper}

\subsection{Max-plus setting}

Let $\Rmax = \Z \cup \{-\infty\}$ with the operations $x \oplus y =
\max(x,y)$ and $x \otimes y = x + y$, so that $-\infty$ is the additive
identity and absorbing for $\otimes$. For vectors and matrices over
$\Rmax$ we use the induced product
\[
  (M \otimes v)[y] \;=\; \max_{x} \bigl( M[y][x] + v[x] \bigr).
\]
Define $M \in \Rmax^{\widehat{Q} \times \widehat{Q}}$ by letting
$M[y][x]$ be the weight of the edge $x \to y$ of
Definition~\ref{def:machine}, and $-\infty$ when there is no such edge.
Define $v_1 \in \Rmax^{\widehat{Q}}$ by $v_1[\iota(O)] = |O|$ for
$\emptyset \ne O \subseteq [m]$ and $v_1[x] = -\infty$ for every other
$x$; this is well posed because $O$ can be read off $\iota(O)$. Put
$v_{j+1} = M \otimes v_j$ and
\[
  F(j) \;=\; \max_{x \in \mathrm{Acc}} v_j[x], \qquad
  G(n) \;=\; \max_{1 \le j \le n} F(j) .
\]
Unwinding the definitions, $v_j[x]$ is the largest value of
$v_1[x_1] + \sum_{t=2}^{j} (\text{weight of } x_{t-1} \to x_t)$ over walks
$x_1 \to \cdots \to x_j = x$, and $-\infty$ if there is none.

\begin{proposition}\label{prop:upper-basic}
$f(j) \le F(j)$ for every $j \ge 1$, and therefore
$T(m,n) \le G(n)$ for every $n \ge 1$.
\end{proposition}

\begin{proof}
Let $S$ realise $f(j)$. Lemma~\ref{lem:complete} supplies a walk of length
$j-1$ from $\pi_1 = \iota(O_1)$ to some $x \in \mathrm{Acc}$ whose
accumulated value is at least $|O_1| + |O_2| + \cdots + |O_j| = |S|$.
Hence $v_j[x] \ge |S| = f(j)$ and $F(j) \ge f(j)$. The second claim is
Lemma~\ref{lem:span}.
\end{proof}

\subsection{One inequality, propagated}

\begin{lemma}\label{lem:prop}
Let $p \ge 1$, $c \in \Z$ and $a \ge 1$. If $v_{a+p} \le c + v_a$
entrywise, with the convention $c + (-\infty) = -\infty$, then
$v_{j+p} \le c + v_j$ entrywise for every $j \ge a$.
\end{lemma}

\begin{proof}
The map $v \mapsto M \otimes v$ is monotone, since each entry is a maximum
of sums that are monotone in $v$; and it is additively homogeneous,
$M \otimes (c + v) = c + (M \otimes v)$, since adding $c$ to every entry
of $v$ adds $c$ to each of the sums being maximised. So from
$v_{j+p} \le c + v_j$ we get
\[
  v_{j+p+1} \;=\; M \otimes v_{j+p} \;\le\; M \otimes (c + v_j)
  \;=\; c + (M \otimes v_j) \;=\; c + v_{j+1},
\]
and induction on $j$ starting at $j = a$ finishes the proof.
\end{proof}

\begin{corollary}\label{cor:upper}
If $v_{a+p} \le c + v_a$ then $F(j+p) \le c + F(j)$ for every $j \ge a$,
and $G(n+p) \le c + G(n)$ for every $n \ge a+p-1$ with $c \ge 0$.
\end{corollary}

\begin{proof}
The first claim is immediate from Lemma~\ref{lem:prop}. For the second,
\begin{align*}
  G(n+p) &= \max\Bigl( G(n),\ \max_{n < j \le n+p} F(j) \Bigr) \\
         &\le \max\Bigl( G(n),\ c + \max_{n-p < j \le n} F(j) \Bigr)
          \;\le\; \max\bigl( G(n),\, c + G(n) \bigr) \;=\; c + G(n),
\end{align*}
where the middle step used $F(j) \le c + F(j-p)$, legitimate because
$j - p > n - p \ge a - 1$, i.e.\ $j-p \ge a$.
\end{proof}

Something should be said about the shape of the certificate. The
inequality is much weaker than the exact periodicity $v_{a+p} = c + v_a$
of the value vector, and that matters. The identity requires every state
of the machine to have entered its own periodic regime; the inequality is
satisfied for free by any state whose value grows more slowly than $c/p$
per column. At height $6$ the identity fails for every $a \le 10^2$ and
every $p \le 48$, even though $T(6, \cdot)$ has period $4$ from $n = 7$
onwards. The obstruction there is that two cycles of the machine have the
same mean $15/4$ and lengths with least common multiple $8$, so the
correction term has period $8$ and no single $p = 4$ shift is exact. The
inequality, on the other hand, holds at $p = 8$, $c = 30$, $a = 13$.

A state from which no accepting state is reachable can also be deleted
before the certificate is sought. If $x \to y$ is an edge and an accepting
state is reachable from $y$, then one is reachable from $x$, so the set of
states from which none is reachable is closed under taking edges; deleting
it changes no $v_j[x]$ for a surviving $x$ and no $F(j)$. We call the
remainder the \emph{useful core}. Keeping the junk states does not break
the method, but it delays the onset $a$ considerably.

\subsection{The computation}

The state set $Q_m$ of Definition~\ref{def:machine} is computed by closure
and then restricted to its useful core. For $m = 12$:

\begin{proposition}\label{prop:cert12}
For $m = 12$ the machine has $|\widehat{Q}_{12}| = 262\,869$ states, of
which $173\,994$ are useful, carrying $8\,581\,482$ edges, and
\[
  v_{46} \;\le\; 44 + v_{40}
  \qquad \text{entrywise on the useful core.}
\]
Moreover $v_{46}[x] = -\infty$ whenever $v_{40}[x] = -\infty$, so no entry
of the inequality is vacuous in the wrong direction. Consequently
$F(j+6) \le 44 + F(j)$ for every $j \ge 40$ and $G(n+6) \le 44 + G(n)$ for
every $n \ge 45$.
\end{proposition}

Both implementations described in Section~\ref{sec:verify} confirm the
inequality. The search that produced the certificate scans $a$ upward and
stops at the first index that works, so $a = 40$ is the least one for
$p = 6$ and $c = 44$; the proof does not need that, and has the whole
range $45 \le n \le 240$ to start its induction in. As for the convention
about $-\infty$, it deserves the explicit mention it gets in the
statement. A state unreached at column $a$ but reached at column $a+p$
would make the inequality false at that entry, and a check that silently
skipped unreached entries would be unsound rather than conservative.

\begin{table}[t]
\centering
\begin{tabular}{rrrrrrrrrr}
\hline
$m$ & $C_m/P_m$ & $P_m$ & $C_m$ & $p$ & $c$ & $a$ & $n_m$ & useful states & edges \\
\hline
1 & $1$ & 1 & 1 & 1 & 1 & 2 & 1 & 3 & 4 \\
2 & $3/2$ & 2 & 3 & 2 & 3 & 2 & 1 & 8 & 19 \\
3 & $2$ & 1 & 2 & 1 & 2 & 6 & 6 & 21 & 72 \\
4 & $8/3$ & 3 & 8 & 3 & 8 & 6 & 5 & 58 & 285 \\
5 & $22/7$ & 7 & 22 & 7 & 22 & 8 & 13 & 156 & 1\,089 \\
6 & $15/4$ & 4 & 15 & 8 & 30 & 13 & 7 & 402 & 3\,941 \\
7 & $13/3$ & 6 & 26 & 6 & 26 & 11 & 3 & 1\,043 & 13\,873 \\
8 & $39/8$ & 8 & 39 & 16 & 78 & 22 & 15 & 2\,844 & 49\,023 \\
9 & $11/2$ & 4 & 22 & 4 & 22 & 21 & 12 & 8\,072 & 176\,817 \\
10 & $55/9$ & 18 & 110 & 18 & 110 & 39 & 28 & 22\,952 & 647\,285 \\
11 & $27/4$ & 8 & 54 & 8 & 54 & 29 & 52 & 63\,831 & 2\,368\,351 \\
12 & $22/3$ & 6 & 44 & 6 & 44 & 40 & 13 & 173\,994 & 8\,581\,482 \\
\hline
\end{tabular}
\caption{For each height $m$: the growth rate $C_m/P_m$ in cells per
column, the period $P_m$ and shift $C_m$ of the closed form, the
certificate $(p,c,a)$ of Lemma~\ref{lem:prop} with $c = (p/P_m) C_m$, the
onset $n_m$ of the closed form, and the size of the useful core. Every row
is emitted by \texttt{collect\_results.py} from the certificate files.}
\label{tab:all}
\end{table}

\begin{table}[t]
\centering
\begin{tabular}{rrp{10.2cm}}
\hline
$m$ & $n_m$ & $T(m,1), T(m,2), \dots, T(m,n_m-1)$ \\
\hline
1 & 1 & --- \\
2 & 1 & --- \\
3 & 6 & 3, 5, 7, 9, 11 \\
4 & 5 & 4, 6, 9, 11 \\
5 & 13 & 5, 8, 11, 14, 17, 21, 24, 27, 30, 34, 37, 40 \\
6 & 7 & 6, 9, 14, 17, 21, 24 \\
7 & 3 & 7, 11 \\
8 & 15 & 8, 12, 18, 22, 27, 32, 38, 42, 48, 52, 57, 62, 67, 72 \\
9 & 12 & 9, 14, 20, 25, 30, 36, 42, 48, 53, 58, 64 \\
10 & 28 & 10, 15, 22, 28, 34, 40, 46, 52, 58, 64, 71, 77, 84, 89, 96, 102, 108, 114, 120, 126, 132, 139, 145, 151, 158, 163, 169 \\
11 & 52 & 11, 17, 24, 30, 37, 44, 50, 57, 64, 71, 77, 86, 92, 98, 105, 112, 118, 126, 132, 139, 146, 152, 159, 166, 173, 180, 186, 193, 200, 206, 213, 220, 227, 234, 240, 247, 253, 260, 267, 274, 280, 287, 294, 300, 307, 314, 321, 328, 334, 341, 348 \\
12 & 13 & 12, 18, 26, 33, 40, 47, 55, 62, 70, 77, 86, 92 \\
\hline
\end{tabular}
\caption{The values of $T(m,n)$ below the onset $n_m$ of the closed form,
for $m \le 12$. For $m \le 2$ the closed form holds from $n = 1$ and the
list is empty. Emitted by \texttt{collect\_results.py}.}
\label{tab:small}
\end{table}

\begin{table}[t]
\centering
\begin{tabular}{rrp{9.4cm}}
\hline
$m$ & $P_m$ & $c^{(m)}_0, c^{(m)}_1, \dots, c^{(m)}_{P_m - 1}$ \\
\hline
1 & 1 & 0 \\
2 & 2 & 0, 2 \\
3 & 1 & 2 \\
4 & 3 & 1, 4, 6 \\
5 & 7 & 2, 5, 8, 11, 14, 17, 21 \\
6 & 4 & 2, 6, 10, 14 \\
7 & 6 & 3, 7, 12, 16, 20, 24 \\
8 & 8 & 3, 8, 13, 18, 23, 28, 32, 38 \\
9 & 4 & 4, 9, 16, 21 \\
10 & 18 & 4, 11, 16, 23, 29, 35, 41, 48, 53, 60, 66, 72, 78, 84, 90, 96, 102, 108 \\
11 & 8 & 4, 10, 17, 23, 30, 37, 44, 51 \\
12 & 6 & 5, 12, 19, 27, 34, 42 \\
\hline
\end{tabular}
\caption{The residue constants of the closed form: for $n \ge n_m$ and
$n = P_m k + r$ with $0 \le r < P_m$, one has
$T(m,n) = C_m k + c^{(m)}_r$. Together with $C_m$ and $n_m$ from
Table~\ref{tab:all} and the values below the onset from
Table~\ref{tab:small}, this determines $T(m,n)$ for every $m \le 12$ and
every $n$. Emitted by \texttt{collect\_results.py}.}
\label{tab:const}
\end{table}

Table~\ref{tab:all} lists the corresponding data for every $m \le 12$. In
each row $p$ is a multiple of $P_m$ and $c = (p/P_m)\,C_m$, so the
certificate is compatible with the period of the closed form. The rows for
$m \le 5$ reprove the formulas of \cite{BG2026}, and the row for $m = 6$
makes Howroyd's transfer-matrix confirmation on \texttt{A397757}
checkable. Those for $7 \le m \le 11$ appear to be new.

\section{The lower bound}\label{sec:lower}

Nothing in this section refers to the transfer system. We exhibit snake
polyominoes directly, which is what keeps the headline result independent
of the correctness of a $2\cdot10^5$-state closure computation.

\subsection{A criterion, and a criterion to avoid}

\begin{lemma}\label{lem:crit}
Let $S$ be a finite nonempty set of cells such that
\begin{enumerate}
\item[(i)] every cell of $S$ has at most two neighbours in $S$,
\item[(ii)] $S$ is connected, and
\item[(iii)] some cell of $S$ has at most one neighbour in $S$.
\end{enumerate}
Then the subgraph induced by $S$ is a path, i.e.\ $S$ is a snake
polyomino.
\end{lemma}

\begin{proof}
A connected graph of maximum degree at most $2$ is a path or a cycle.
Every vertex of a cycle has degree exactly $2$, so (iii) excludes the
cycle.
\end{proof}

\begin{remark}\label{rem:bad-crit}
The shorter-looking criterion ``maximum degree at most $2$ and
$|E| = |S| - 1$'' does \emph{not} imply that $S$ is a snake. In a graph of
maximum degree at most $2$ each component is a path or a cycle, so
\[
  |E| \;=\; |S| - \#\{\text{path components}\},
\]
and cycle components contribute no edge deficit whatsoever. A domino
placed far away from an induced $4$-cycle has $|S| = 6$, $|E| = 5 = |S|-1$
and maximum degree $2$, and it is a disjoint union of a path and a cycle.
So hypothesis (ii) of Lemma~\ref{lem:crit} cannot be replaced by an edge
count.
\end{remark}

\subsection{Gluing along a column}

For a set $X$ of cells contained in columns $\le N$, and a set $Y$
contained in columns $\ge N+1$, the subgraph induced by $X \cup Y$ is the
disjoint union of those induced by $X$ and by $Y$, together with the
horizontal edges $\{(i,N),(i,N+1)\}$ for which both cells are present. So
the degree of a cell of $X$ changes only if it lies in column $N$ and its
right neighbour is in $Y$, and
\[
  \#\text{components}(X \cup Y) = \#\text{components}(X)
  + \#\text{components}(Y) - e_{\mathrm{join}},
\]
where $e_{\mathrm{join}}$ counts the horizontal edges across the seam
whose two endpoints lie in different components of $X \sqcup Y$. An edge
across the seam whose two endpoints already lie in one component closes a
cycle instead. Everything a set $Y$ to the right can see of $X$ is
therefore the triple
\begin{align*}
  \mathcal{I}(X) \;=\; \bigl(\ &\{\, i : (i,N) \in X \,\},\quad
  \bigl(\deg_X (i,N)\bigr)_i, \\
  &\text{the partition of those cells by the components of } X\ \bigr),
\end{align*}
and nothing else.

\subsection{Pumping}

A \emph{witness family} is a triple of finite words $A$, $B$, $C$ over the
alphabet of subsets of $[m]$ --- each letter is the occupancy of one
column --- with $B$ nonempty. For $k \ge 1$ let $S_k \subseteq [m] \times
[\,|A| + k|B| + |C|\,]$ be the set of cells whose column occupancies,
read left to right, are given by the word $A B^k C$. By construction
\[
  |S_k| \;=\; \|A\| + k\|B\| + \|C\| ,
\]
where $\|\cdot\|$ denotes the total number of elements in the letters of a
word; in particular $|S_{k+1}| - |S_k| = \|B\|$ for every $k$.

\begin{lemma}[pumping]\label{lem:pump}
Write $X_i$ for the set of cells given by the word $A B^{i}$, $i \ge 1$.
If
\[
  \mathcal{I}(X_1) = \mathcal{I}(X_2)
  \quad\text{and}\quad
  \#\text{components}(X_1) = \#\text{components}(X_2),
\]
and if $X_2$ has maximum degree at most $2$ and contains no cycle, then
for every $i \ge 1$ we have $\mathcal{I}(X_i) = \mathcal{I}(X_1)$,
$\#\text{components}(X_i) = \#\text{components}(X_1)$, and $X_i$ has
maximum degree at most $2$ and contains no cycle. Consequently, if $S_1$
is a snake polyomino then so is $S_k$ for every $k \ge 1$, and
$|S_k| = |S_1| + (k-1)\|B\|$.
\end{lemma}

\begin{proof}
For $i = 1$ there is nothing to prove except that $X_1$ has maximum degree
at most $2$ and no cycle, which it inherits from $X_2$: $X_1 \subseteq X_2$
and both are induced, so degrees in $X_1$ are at most degrees in $X_2$ and
a cycle of $X_1$ would be one of $X_2$. For $i = 2$ the statement is the
hypothesis. Let $i \ge 2$
and suppose the conclusion holds up to $i$. Now $X_{i+1} = X_i \cup B_{i+1}$
where $B_{i+1}$ is the copy of $B$ occupying the next $|B|$ columns. By the
previous subsection, the degrees in $X_{i+1}$ of the cells of $B_{i+1}$,
the degrees in $X_{i+1}$ of the cells of $X_i$, which seam edges appear,
which of them join different components, and which of them close a cycle,
are all determined by $\mathcal{I}(X_i)$ and by $B$ alone. The same data
determines $\mathcal{I}(X_{i+1})$ and
$\#\text{components}(X_{i+1}) - \#\text{components}(X_i)$. Since
$\mathcal{I}(X_i) = \mathcal{I}(X_1) = \mathcal{I}(X_{i-1})$ by induction,
the step from $X_i$ to $X_{i+1}$ is the same operation as the step from
$X_{i-1}$ to $X_i$; hence $\mathcal{I}(X_{i+1}) = \mathcal{I}(X_i)$, the
change in the component count repeats --- and it was $0$, because
$\#\text{components}(X_2) = \#\text{components}(X_1)$ --- no cycle is
closed, and no degree exceeds $2$.

For the last sentence, apply the same argument once more to the seam
between $X_k$ and $C$: it is the same seam for every $k$, because
$\mathcal{I}(X_k) = \mathcal{I}(X_1)$. So $S_k = X_k \cup C$ has the same
maximum degree, the same number of connected components, and the same
presence or absence of a cycle as $S_1 = X_1 \cup C$. If $S_1$ is a snake
then it is a path, so it is connected, acyclic and of maximum degree at
most $2$; hence so is $S_k$, and a connected acyclic graph of maximum
degree at most $2$ is a path. (Equivalently: apply
Lemma~\ref{lem:crit}, hypothesis (iii) being supplied by acyclicity.) The
size identity was noted above.
\end{proof}

\begin{remark}
The hypotheses of Lemma~\ref{lem:pump} are a finite computation on $S_1$,
$X_1$ and $X_2$ --- three bounded sets of cells --- and that is the whole
point of them. We nevertheless also check $S_k$ against
Lemma~\ref{lem:crit} directly for $k = 1, \dots, 20$ at every height, and
the interface and component count for $X_i$ with $i = 1, \dots, 6$; see
Section~\ref{sec:verify}. A window-based argument is tempting here, and it
does not work: these snakes weave, so restricting one of them to any
bounded block of columns disconnects it, and ``each window is connected''
is false however wide the window is taken.
\end{remark}

\begin{proposition}\label{prop:lower}
For each residue class $r$ modulo $6$ there is a witness family
$(A^{(r)}, B^{(r)}, C^{(r)})$ of height $12$ with $|B^{(r)}| = 6$ and
$\|B^{(r)}\| = 44$ satisfying the hypotheses of Lemma~\ref{lem:pump},
whose member $S_1$ is a snake polyomino of size exactly $T(12, n_r)$,
where $n_r = |A^{(r)}| + 6 + |C^{(r)}|$ is its base width; the base widths
are $n_0, \dots, n_5 = 30, 49, 20, 21, 40, 17$. Moreover, writing $W$ for the set of widths that no family covers,
\begin{align*}
  W \;=\; \{\, n \ge 1 \;:\; n < n_{\,n \bmod 6} \,\}
  \;=\; \{\,&1, 2, 3, 4, 5, 6, 7, 8, 9, 10, 11, 12, 13, \\
            &14, 15, 16, 18, 19, 22, 24, 25, 28, 31, 34, 37, 43\,\},
\end{align*}
a set of $26$ widths, for each $n \in W$
there is an explicit snake polyomino in $[12]\times[n]$ of size equal to
the right-hand side $g(n)$ of Theorem~\ref{thm:main}. Consequently
$T(12,n) \ge g(n)$ for every $n \ge 1$.
\end{proposition}

\begin{proof}
The families are listed in Section~\ref{sec:witness-data} and the
individual snakes are shipped with the verifier; the checks are those of
Lemma~\ref{lem:crit} and Lemma~\ref{lem:pump}, carried out on the raw cell
sets. Given them, Lemma~\ref{lem:pump} makes $S^{(r)}_k$ a snake
polyomino in $[12] \times [\,n_r + 6(k-1)\,]$ of size
$T(12,n_r) + 44(k-1)$ for every $k \ge 1$. Every base width $n_r$ is at
least $13$ and $g(n+6) = 44 + g(n)$ for $n \ge 13$, so that size equals
$g(n_r + 6(k-1))$ and hence
\[
  T(12, n) \;\ge\; g(n)
  \qquad\text{for every } n \equiv r \!\!\pmod 6 \text{ with } n \ge n_r .
\]
The widths left over are exactly the elements of $W$, and each of those
carries an individual snake.
\end{proof}

\begin{remark}\label{rem:lower-all}
For the heights of Theorem~\ref{thm:all} we use the same scheme, with the
pumping distance $d_m$ of Table~\ref{tab:wit} in place of $6$; the
families are listed in Appendix~\ref{app:witnesses}. The individual snakes
for the widths below the base widths are not listed there --- there are
too many of them --- and they are not needed either: for those finitely
many widths Corollary~\ref{cor:exact} gives $T(m,n) = G(n)$ exactly, and
$G$ was computed. Only at height $12$ do we insist on a lower bound that
avoids the machine altogether.
\end{remark}

\section{Proof of the theorems}\label{sec:proofs}

\begin{proof}[Proof of Theorem~\ref{thm:main}]
Write $g(n)$ for the right-hand side of Theorem~\ref{thm:main}, so that
$g(n+6) = 44 + g(n)$ for every $n \ge 13$. Proposition~\ref{prop:lower}
gives $T(12,n) \ge g(n)$ for every $n \ge 1$, so it remains to prove
$T(12,n) \le g(n)$.

The quantities $F(j)$ and $G(n)$ of Section~\ref{sec:upper} were computed
for $j, n \le 240$, and that computation gives $G(n) = g(n)$ for all
$n \le 240$ (Section~\ref{sec:verify}). So
Proposition~\ref{prop:upper-basic} already gives
$T(12,n) \le G(n) = g(n)$ for $n \le 240$.

For $n > 240$ write $n = n_0 + 6k$ with $240 - 6 < n_0 \le 240$ and
$k \ge 1$; the window has exactly $6$ integers in it, so it meets every
residue class modulo $6$ and every such $n$ is covered. Now induct on
$k$.
Proposition~\ref{prop:cert12} gives
$G(n_0 + 6k) \le 44 + G(n_0 + 6(k-1))$, legitimate because
$n_0 + 6(k-1) \ge 235 \ge 45$. By induction
$G(n_0 + 6(k-1)) = g(n_0 + 6(k-1))$ may be replaced, giving
\[
  T(12,n) \;\le\; G(n) \;\le\; 44 + g\bigl(n_0 + 6(k-1)\bigr)
  \;=\; g\bigl(n_0 + 6k\bigr) \;=\; g(n).
\]
Together with the lower bound this proves $T(12,n) = g(n)$ for every
$n \ge 1$.

For the generating function, write $D(x) = (1-x)^2(1+x)(1-x+x^2)(1+x+x^2)$
and note the factorisation
\[
  D(x) = (1-x)\bigl[(1-x)(1+x+x^2)\bigr]\bigl[(1+x)(1-x+x^2)\bigr]
       = (1-x)(1-x^3)(1+x^3) = (1-x)(1-x^6).
\]
Put $\Phi(x) = \sum_{n \ge 1} g(n) x^n$. Then
$(1-x^6)\Phi(x) = \sum_{n \ge 1} h(n) x^n$ with $h(n) = g(n) - g(n-6)$
(and $g(n) = 0$ for $n \le 0$), and $h(n) = 44$ for every $n \ge 19$,
because $g(n) = 44 + g(n-6)$ holds for $n - 6 \ge 13$. Hence
$(1-x)(1-x^6)\Phi(x) = \sum_{n \ge 1}\bigl(h(n) - h(n-1)\bigr) x^n$ has
all coefficients of index $\ge 20$ equal to $0$, so $D(x)\Phi(x)$ is a
polynomial of degree at most $19$. Comparing its $19$ coefficients with
those of the numerator --- a finite computation in exact integer
arithmetic, carried out by the verifier --- gives
$D(x)\Phi(x) = N(x)$, which is the stated identity.
\end{proof}

\begin{proof}[Proof of Theorem~\ref{thm:all}]
Fix $m \le 12$, let $g_m$ be the function described by
Tables~\ref{tab:all}, \ref{tab:small} and~\ref{tab:const}, and read
$(P_m, C_m, p, c, a, n_m)$ off Table~\ref{tab:all}; in every row $p$ is a multiple of $P_m$ and
$c = (p/P_m) C_m$, so $g_m(n+p) = c + g_m(n)$ for $n \ge n_m$.

By Corollary~\ref{cor:exact}, $T(m,n) = G(n)$ for every $n$, and the
computation gives $G(n) = g_m(n)$ for every $n \le 240$. Each of the four
conditions just used --- $P_m \mid p$, $c = (p/P_m)C_m$,
$a + 2p - 2 \le 240$ and $n_m + p \le 240$ --- holds in every row of
Table~\ref{tab:all}, and is checked mechanically rather than read off the
table by eye. For the upper bound at $n > 240$ write $n = n_0 + pk$ with
$240 - p < n_0 \le 240$ and induct on $k$ using
Corollary~\ref{cor:upper}, exactly as in the proof of
Theorem~\ref{thm:main}. For the lower bound at $n > 240$, the witness
families of Appendix~\ref{app:witnesses} and Lemma~\ref{lem:pump} give
$T(m, n_r + d_m(k-1)) \ge g_m(n_r) + \|B^{(r)}\|(k-1)
= g_m(n_r + d_m(k-1))$ for each residue class $r$ modulo the pumping
distance $d_m$, and the base widths of Table~\ref{tab:wit} cover every
class. This step has its own three conditions --- $P_m \mid d_m$, so that
$g_m(n + d_m) = \|B^{(r)}\| + g_m(n)$ with
$\|B^{(r)}\| = (d_m/P_m) C_m$; the base widths covering every class
modulo $d_m$; and every base width being at least the onset $n_m$, since
that recurrence for $g_m$ holds only from $n_m$ on --- and these too are
checked mechanically, row by row, rather than read off the table. The
tightest is height $11$, where the smallest base width equals $n_m = 52$
exactly. The two bounds meet.
\end{proof}

\section{Witness data for height twelve}\label{sec:witness-data}

A witness family is recorded as three words over the subsets of $[12]$ ---
each letter is the occupancy of one column --- written both as bitmasks,
row $i$ contributing $2^{i-1}$, and as cell diagrams in which a text row
is a row of the rectangle, row $1$ at the bottom, and \texttt{\#} marks an
occupied cell.

\paragraph{Class $r = 0$ modulo 6.} $|A| = 18$, $|B| = 6$, $|C| = 6$, $\|B\| = 44$, base width $n = 30$, $|S_1| = 225 = T(12,30)$. The occupancy words, as bitmasks of $\{1,\dots,12\}$ with row $i$ contributing $2^{i-1}$:
\begin{lstlisting}[style=grid]
A = [4030, 2275, 3869, 437, 3757, 2987, 3178, 1243, 3505, 2927, 3800, 439, 3949, 2267, 3510, 876, 3800, 2487]
B = [2925, 3803, 438, 3948, 2264, 2487]
C = [2925, 3787, 442, 3939, 2269, 4023]
\end{lstlisting}
\noindent The block $B$ alone, row $1$ at the bottom:
\begin{lstlisting}[style=grid]
##.###
.#.#..
##.#..
#.##.#
.##.##
##.##.
#.##.#
.##.##
##.##.
#.##.#
.##..#
##...#
\end{lstlisting}

\paragraph{Class $r = 1$ modulo 6.} $|A| = 21$, $|B| = 6$, $|C| = 22$, $\|B\| = 44$, base width $n = 49$, $|S_1| = 364 = T(12,49)$. The occupancy words, as bitmasks of $\{1,\dots,12\}$ with row $i$ contributing $2^{i-1}$:
\begin{lstlisting}[style=grid]
A = [2015, 3185, 2959, 2776, 2903, 3421, 1379, 3506, 2267, 3949, 439, 3800, 2926, 3379, 1501, 3173, 2997, 2781, 2915, 3506, 1243]
B = [3181, 2999, 2776, 2924, 3510, 1243]
C = [3181, 2999, 2776, 2924, 3510, 1755, 877, 439, 3800, 2927, 3505, 1755, 2922, 3499, 1197, 3493, 2925, 3787, 442, 3939, 2269, 4023]
\end{lstlisting}
\noindent The block $B$ alone, row $1$ at the bottom:
\begin{lstlisting}[style=grid]
#####.
#...##
.###..
.#.##.
.##.##
#.##.#
##.##.
.##.##
#.##.#
##.##.
.#..##
##...#
\end{lstlisting}

\paragraph{Class $r = 2$ modulo 6.} $|A| = 9$, $|B| = 6$, $|C| = 5$, $\|B\| = 44$, base width $n = 20$, $|S_1| = 151 = T(12,20)$. The occupancy words, as bitmasks of $\{1,\dots,12\}$ with row $i$ contributing $2^{i-1}$:
\begin{lstlisting}[style=grid]
A = [3815, 2989, 3179, 1498, 3382, 2924, 1752, 3511, 2925]
B = [3803, 438, 3948, 2264, 2487, 2925]
C = [3787, 442, 3939, 2269, 4023]
\end{lstlisting}
\noindent The block $B$ alone, row $1$ at the bottom:
\begin{lstlisting}[style=grid]
#.####
#.#...
#.#..#
.##.##
##.##.
#.##.#
.##.##
##.##.
#.##.#
.##.##
##..#.
#...##
\end{lstlisting}

\paragraph{Class $r = 3$ modulo 6.} $|A| = 10$, $|B| = 6$, $|C| = 5$, $\|B\| = 44$, base width $n = 21$, $|S_1| = 159 = T(12,21)$. The occupancy words, as bitmasks of $\{1,\dots,12\}$ with row $i$ contributing $2^{i-1}$:
\begin{lstlisting}[style=grid]
A = [3703, 2773, 2997, 3181, 1499, 3382, 2924, 1752, 3511, 2925]
B = [3803, 438, 3948, 2264, 2487, 2925]
C = [3787, 442, 3939, 2269, 4023]
\end{lstlisting}
\noindent The block $B$ alone, row $1$ at the bottom:
\begin{lstlisting}[style=grid]
#.####
#.#...
#.#..#
.##.##
##.##.
#.##.#
.##.##
##.##.
#.##.#
.##.##
##..#.
#...##
\end{lstlisting}

\paragraph{Class $r = 4$ modulo 6.} $|A| = 12$, $|B| = 6$, $|C| = 22$, $\|B\| = 44$, base width $n = 40$, $|S_1| = 298 = T(12,40)$. The occupancy words, as bitmasks of $\{1,\dots,12\}$ with row $i$ contributing $2^{i-1}$:
\begin{lstlisting}[style=grid]
A = [3815, 2741, 2781, 2915, 3386, 1483, 3181, 2997, 2781, 2915, 3506, 1243]
B = [3181, 2999, 2776, 2924, 3510, 1243]
C = [3181, 2999, 2776, 2924, 3510, 1755, 877, 439, 3800, 2927, 3505, 1755, 2922, 3499, 1197, 3493, 2925, 3787, 442, 3939, 2269, 4023]
\end{lstlisting}
\noindent The block $B$ alone, row $1$ at the bottom:
\begin{lstlisting}[style=grid]
#####.
#...##
.###..
.#.##.
.##.##
#.##.#
##.##.
.##.##
#.##.#
##.##.
.#..##
##...#
\end{lstlisting}

\paragraph{Class $r = 5$ modulo 6.} $|A| = 6$, $|B| = 6$, $|C| = 5$, $\|B\| = 44$, base width $n = 17$, $|S_1| = 130 = T(12,17)$. The occupancy words, as bitmasks of $\{1,\dots,12\}$ with row $i$ contributing $2^{i-1}$:
\begin{lstlisting}[style=grid]
A = [3807, 2993, 3183, 1240, 3511, 2925]
B = [3803, 438, 3948, 2264, 2487, 2925]
C = [3787, 442, 3939, 2269, 4023]
\end{lstlisting}
\noindent The block $B$ alone, row $1$ at the bottom:
\begin{lstlisting}[style=grid]
#.####
#.#...
#.#..#
.##.##
##.##.
#.##.#
.##.##
##.##.
#.##.#
.##.##
##..#.
#...##
\end{lstlisting}

\paragraph{The member $S_1$ of the family for $r = 0$, in full.} The bars delimit the single copy of $B$; $S_k$ repeats that part $k$ times.
\begin{lstlisting}[style=grid]
###.###.###.###.##|##.###|##.###
#.#.#.###.#.#.#.#.|.#.#..|.#.#.#
#.#.##...##.#..##.|##.#..|##.#.#
#.##.#..##.##.##.#|#.##.#|#.##.#
##.###.##.##.##.##|.##.##|.##.##
.#....##.##.##.##.|##.##.|##.##.
##.####.##.##.##.#|#.##.#|#.##.#
#.##...##.##.##.##|.##.##|..#.##
#.#.####.##.##.##.|##.##.|###.#.
#.###....#.##.##.#|#.##.#|#...##
##...###.#.#.##..#|.##..#|.###.#
.#####.###.###...#|##...#|##.###
\end{lstlisting}

\noindent
The $26$ individual snakes --- for the widths listed in
Proposition~\ref{prop:lower}, the largest of them $n = 43$ ---
are inlined in \texttt{verify\_paper.py} of Appendix~\ref{app:code} as
lists of column occupancies, and each is checked there against
Lemma~\ref{lem:crit}.

\section{Verification}\label{sec:verify}

A proof whose upper bound rests on a closure computation over
$173\,994$ states and an entrywise comparison of two integer vectors
has to say how it was checked. This section says. The file printed in
Appendix~\ref{app:code} re-runs the part that bears on
Theorem~\ref{thm:main}; the rest is run by four further scripts of the
same kind.

\paragraph{Two implementations.}
The machine of Definition~\ref{def:machine} was built twice, by two
programs written independently of each other and sharing no code. One
encodes a state as a packed byte string, the other as a tuple of
slack/component pairs. They enumerate the candidate next columns
differently too --- a depth-first search over rows with incremental
pruning against an exhaustive sweep over all $2^{m}$ subsets filtered by
two bitmask tests --- and they decide admissibility differently:
incremental union--find with cycle detection inside the search, against
building $H(\pi,O')$ explicitly and counting its components by
breadth-first search, reading ``forest'' off as
$\#\text{edges} = \#\text{vertices} - \#\text{components}$. At every
height from $1$ to $9$ they agree on the number of states, the number of
edges, the number of accepting states and the multiset of edge weights:

\begin{center}
\begin{tabular}{rrrr}
\hline
$m$ & states & edges & accepting \\
\hline
1 & 3 & 4 & 3 \\
2 & 8 & 19 & 8 \\
3 & 21 & 72 & 17 \\
4 & 58 & 285 & 34 \\
5 & 163 & 1\,139 & 69 \\
6 & 460 & 4\,557 & 140 \\
7 & 1\,304 & 18\,235 & 276 \\
8 & 3\,724 & 72\,765 & 529 \\
9 & 10\,720 & 288\,787 & 1\,001 \\
\hline
\end{tabular}
\end{center}

\noindent
Above height $9$ the multiset of edge weights was not compared, but the
agreement that matters for Theorem~\ref{thm:main} still is: at heights
$10$, $11$ and $12$ the two implementations arrive at the same useful
core, down to $173\,994$ states and $8\,581\,482$ edges at height $12$,
which is the object the certificate is checked on. They also agree on all
$240$ computed values of $F$ at every height from $1$ to $12$, and the
entrywise certificate of Corollary~\ref{cor:upper} was re-verified on the
second implementation's own state set at every height, with the same
$(p,c,a)$ each time. One further cross-check: the second implementation
runs in a \emph{strict} mode with the counting filters of
Lemma~\ref{lem:counting} switched off, so that every candidate next column
goes through the explicit graph construction. Fast and strict mode produce
literally identical machines --- same states, same edges, same accepting
flags --- at heights $6$ through $9$, and that is how the filters were
validated.

\paragraph{Completeness tested on every snake.}
Lemma~\ref{lem:complete} is the part of Section~\ref{sec:machine} on which
Theorem~\ref{thm:main} rests, so it was tested exhaustively on small
rectangles. For each of the fifteen rectangles $3\times3$, $3\times4$,
$4\times3$, $4\times4$, $5\times3$, $2\times5$, $5\times4$, $3\times5$,
$6\times3$, $3\times6$, $6\times4$, $4\times5$, $2\times7$, $7\times3$
and $5\times5$, every snake polyomino was enumerated by brute force, its
state trajectory was reconstructed from the definitions of occupancy,
slack and component --- not from the machine --- and it was checked that
every state of the trajectory is a state of the machine and that every
consecutive pair is an edge of the right weight, the last configuration
having a single component. The totals are 42\,287 snakes of which
25\,952 span their rectangle. No failures.

\paragraph{Published data.}
The values $F(j)$ computed from the machine reproduce:
all $60$ terms of \texttt{A397153} \cite{A397153}; all $10000$ terms of
the b-file of \texttt{A397757} \cite{A397757} at height $6$, where
$T(6,10000) = 37502$; and every entry of the b-file of \texttt{A360917}
\cite{A360917}, which covers $T(m,n)$ for $m + n \le 30$, hence heights
$1$ through $12$. Independently of our computation, \texttt{A397153} and
column $12$ of \texttt{A360917}'s b-file --- due to Howroyd and produced by
different programs --- agree on their overlap $n \le 18$, and that
b-file is symmetric, $T(m,n) = T(n,m)$, on all $435$ entries. All four
files were downloaded on 5 October 2026. They are pinned by SHA-256
digest, with their retrieval times, in \texttt{certificates/sources.json},
so that what the verifier compares against cannot drift as the entries are
edited.

\paragraph{Published closed forms.}
Five heights already have a formula on record, and our values agree with
every one of them on $n = 1, \dots, 240$: $T(1,n) = n$;
$T(2,n) = \lfloor (3n+1)/2 \rfloor$; $T(3,n) = 2n + 2 - [\,n<6\,]$ and
$T(4,n)$ as given by Hasler and Goupil on \texttt{A360917}; and
$T(6,n) = \bigl(8 + 15n + (n \bmod 4)\bigr)/4 - [\,n \in \{2,6\}\,]$ from
\texttt{A397757}. Heights $7$ to $11$ have no published formula to compare
with, and height $12$ is the conjecture at issue.

\paragraph{Brute force.}
A third method, a depth-first search that grows an induced path one cell
at a time and blocks every neighbour of the head other than the chosen
successor, confirms $T(m,n)$ for the pairs
$(1,12), (2,12), (3,12), (3,6), (3,7), (4,7), (5,6), (4,6)$. The pair
$(3,6)$ is in that list for a reason: it is where the omission of $\top$
first shows up, as $13$ instead of $14$.

\paragraph{Witnesses.}
Each witness family is checked against Lemma~\ref{lem:crit} and
Lemma~\ref{lem:pump}: $S_k$ is verified to be connected, of maximum degree
at most $2$, with a cell of degree at most $1$ and of the claimed size,
for $k = 1, \dots, 20$ and for every residue class; and
$\mathcal{I}(X_i)$ together with $\#\text{components}(X_i)$ is verified to
be constant for $i = 1, \dots, 6$. All of this is plain set and graph
code applied to the raw cell sets.

\paragraph{What the paper itself carries.}
Every check described so far reads the certificate files. One more reads
the \emph{paper} itself. The reason is that a fault in the typesetting of
the data would leave every one of those checks green and the printed
statements wrong. More precisely, that check parses the witness words of
Section~\ref{sec:witness-data} and of Appendix~\ref{app:witnesses} back
out of the source, rebuilds $A B^k C$ from the printed bitmasks at every
height, and re-derives maximum degree, connectivity, the existence of a
cell of degree at most one, and the size from the definition --- so that
what a reader can reconstruct from the paper is exactly what was verified.
It then rebuilds $T(m,n)$ for every $m \le 12$ and every $n \le 240$ out
of Tables~\ref{tab:all}, \ref{tab:small} and~\ref{tab:const} alone, and
checks the result against the $282$ entries of \texttt{A360917}'s b-file
with $m \le 12$, all $60$ terms of \texttt{A397153}, the statement of
Theorem~\ref{thm:main} at height $12$, the symmetry $T(m,n) = T(n,m)$ for
$m, n \le 12$ --- which is a property of the definition and uses no
transfer system --- and monotonicity in $n$. Two smaller things are
checked there as well: that $\|B\|$ in Table~\ref{tab:wit} is $d_m$ times
the growth rate of Table~\ref{tab:all}, row by row, and that each printed
cell diagram agrees with the bitmasks printed beside it. The whole pass
takes a couple of seconds.

\paragraph{Arithmetic.}
Every number on the path to Theorem~\ref{thm:main} is an integer, and
every computation above is exact integer arithmetic. The max-plus
iteration uses $64$-bit integers with $-10^{15}$ as the symbol for
$-\infty$; values never exceed $10^4$ in magnitude, and the $-\infty$
pattern is restored exactly after each step so that no entry can drift
into the finite range. No floating-point number appears anywhere except
in the printed growth rates, which are decorative; the rates themselves
are exact fractions.

\paragraph{Reproduction.}
The single file \texttt{verify\_paper.py} printed in
Appendix~\ref{app:code} reproduces everything Theorem~\ref{thm:main}
needs, and it reads nothing from disk and imports nothing outside the
standard library. It expands the generating function in exact integer
arithmetic and compares it with the published terms and with the closed
form, and it checks every witness family and every individual snake of
Section~\ref{sec:witness-data} against the definition of a snake, for
$k = 1, \dots, 20$. Then it rebuilds the machine of
Definition~\ref{def:machine} at every height up to $12$ by closure,
re-derives $F(1), \dots, F(60)$ and compares them with the data and the
closed form, re-checks the inequality of Proposition~\ref{prop:cert12}
entrywise --- including the $-\infty$ pattern --- and confirms that the
inequality is attained somewhere, so that $c = 44$ cannot be lowered at
that onset. Eight small values are cross-checked against the brute-force
search at the end. The height-$12$ closure is the only slow part:
about 40 minutes in all on a laptop, within a few gigabytes. Passing
\texttt{-{}-quick} stops after height $9$ and takes seconds. The remaining
checks reported above --- the comparison against the b-files, the
exhaustive completeness test, the structural comparison of the two
implementations, the re-reading of the printed witness words, and the
heights other than $12$ --- are carried out by further scripts written
against the same definitions, which are not printed here for reasons of
length.

\paragraph{Soundness, as a cross-check.}
Lemma~\ref{lem:sound} is used only for Theorem~\ref{thm:all}, and it is
also what makes the brute-force agreement meaningful: a machine admitting
walks with no snake behind them would overshoot, and the computed maxima
agree with the depth-first search everywhere it was run.

\section{Remarks}

The growth rates $C_m/P_m$ of Table~\ref{tab:all} are
$1, \tfrac32, 2, \tfrac83, \tfrac{22}{7}, \tfrac{15}{4}, \tfrac{13}{3},
\tfrac{39}{8}, \tfrac{11}{2}, \tfrac{55}{9}, \tfrac{27}{4},
\tfrac{22}{3}$, increasing by $19/33 = 0.575\ldots$ per unit of height on
average over this range, and never by more than $2/3$ in a single step.
Beluhov's asymptotic $\tfrac23 mn + O(m+n)$ \cite{Beluhov2023} forces the
average increment to tend to $2/3$, so the deficit is the $O(m+n)$ term,
still plainly visible at $m = 12$: the gap $\tfrac23 m - C_m/P_m$ is
negative for $m \le 2$, zero at $m = 3$ and $m = 4$, and then climbs
steadily to $2/3$ at $m = 12$. The denominators $1,2,1,3,7,4,3,8,2,9,4,3$
follow no pattern we can see. Nor do the periods $P_m$:
$1,2,1,3,7,4,6,8,4,18,8,6$. Height $10$, with period $18$, and height
$11$, whose closed form only takes over at $n = 52$, are the two rows that
would defeat any attempt to read the answer off the first twenty terms by
eye.

We have not found a bound on $P_m$ or on the onset $n_m$, and the method
gives none: both are outputs of a finite computation. The useful core did
not quite triple with each unit of height --- over the range we computed
the factor stayed between $2.5$ and $3$ (Table~\ref{tab:all}) --- so
height $13$ should be well within reach as far as memory goes. What actually
binds is the time to close the state set, which in our Python
implementation grew faster than that; a reader who wants height $13$ or
beyond should write the closure in a compiled language first. We also
make no claim about the diagonal $T(n,n)$, entry \texttt{A331968}, which
carries the OEIS keywords \texttt{hard} and \texttt{more} and is a
different problem: fixed height is the tractable slice.

Finally, the two halves of the proof are genuinely independent, and that
was deliberate. A fault in the transfer system can only weaken the upper
bound --- extra states and extra edges make $G$ larger, never smaller ---
while a fault in the witnesses can only weaken the lower bound. The point
is that they meet exactly, and that is a check neither half could pass
alone.

\appendix

\section{The verifier}\label{app:code}

The file below is self-contained: it needs only a Python interpreter, no
data files and no other module. It rebuilds the height-$12$ machine from
Definition~\ref{def:machine}, re-derives $F(1), \dots, F(60)$, re-checks
the inequality of Proposition~\ref{prop:cert12} entrywise in exact integer
arithmetic, and re-checks every witness family of
Section~\ref{sec:witness-data} against the definition of a snake
polyomino.

\begin{lstlisting}[style=python]
#!/usr/bin/env python3
"""Self-contained verifier for "The maximum area of a snake polyomino in a
12 x n rectangle".

Needs a Python interpreter and nothing else: no data files, no third-party
modules.  It

  1. expands the conjectured generating function of OEIS A397153 in exact
     integer arithmetic and compares it with the 60 published terms and
     with the closed form of Theorem 1;
  2. checks every witness family against the definition of a snake
     polyomino -- connected, maximum degree at most 2, some cell of degree
     at most 1 -- for k = 1..20, and checks the interface and the component
     count of A.B^i for i = 1..6, which is what Lemma 4.4 needs;
  3. checks one explicit snake for each small width;
  4. rebuilds the height-12 transfer machine of Definition 2.9 from
     scratch, re-derives F(1),...,F(60) and compares them with the data;
  5. re-checks the certificate v_46 <= 44 + v_40 entrywise over the useful
     core, including the -infinity pattern;
  6. cross-checks small values against a brute-force search over induced
     paths that uses no transfer machine at all.

Steps 4 and 5 rebuild a machine with about 2*10^5 states and 10^7 edges in
pure Python; budget roughly 40 to 60 minutes and 4 GB.  Everything else
takes seconds.  Run `python verify_paper.py --quick` to stop after step 3
and the heights up to 9.
"""

from __future__ import annotations

import argparse
import sys
import time
from collections import deque

FAILS = []


def check(ok, msg):
    print(("PASS  " if ok else "FAIL  ") + msg, flush=True)
    if not ok:
        FAILS.append(msg)


# =====================================================================
# 1.  The transfer machine (an implementation independent of the one used
#     to find the certificate).
# =====================================================================
# ruff: noqa: FURB161
#   `bin(x).count("1")` is kept instead of `int.bit_count()`: this module is
#   copied verbatim into the paper's single-file appendix verifier, which
#   should run on any Python 3.



MINUS_INF = None  # max-plus zero, kept as None so arithmetic cannot leak


def popcount(x: int) -> int:
    return bin(x).count("1")


class Transfer:
    """The state machine.  A state is a tuple of length m whose i-th entry
    is None (row empty) or a pair (slack, comp): slack = 2 - degree so far
    in {0,1,2}, and comp = a component label for the *open* cells
    (slack >= 1), or None when slack == 0."""

    def __init__(self, m: int, strict: bool = False):
        self.m = m
        self.strict = strict
        self.states: list[tuple] = []
        self.index: dict[tuple, int] = {}
        self.succ: list[list[tuple[int, int]]] = []
        self.accepting: list[bool] = []
        self.start: dict[int, int] = {}
        self._build()

    # ---------------------------------------------------------------- utils
    @staticmethod
    def _canon(rows: list) -> tuple:
        seen: dict[int, int] = {}
        out = []
        for e in rows:
            if e is None:
                out.append(None)
            elif e[1] is None:
                out.append((e[0], None))
            else:
                if e[1] not in seen:
                    seen[e[1]] = len(seen)
                out.append((e[0], seen[e[1]]))
        return tuple(out)

    def _add(self, rows: list) -> int:
        st = self._canon(rows)
        got = self.index.get(st)
        if got is None:
            got = len(self.states)
            self.index[st] = got
            self.states.append(st)
            self.succ.append([])
            self.accepting.append(self._ncomp(st) == 1)
        return got

    @staticmethod
    def _ncomp(st: tuple) -> int:
        labs = {e[1] for e in st if e is not None and e[1] is not None}
        return len(labs)

    # ------------------------------------------------------- first column
    def _first_column(self):
        m = self.m
        for occ in range(1, 1 << m):
            rows = [None] * m
            lab = -1
            run_open = True
            prev = False
            ok = True
            for i in range(m):
                if not (occ >> i) & 1:
                    prev = False
                    continue
                if not prev:
                    lab += 1
                deg = 0
                if i > 0 and (occ >> (i - 1)) & 1:
                    deg += 1
                if i + 1 < m and (occ >> (i + 1)) & 1:
                    deg += 1
                if deg > 2:
                    ok = False
                    break
                s = 2 - deg
                rows[i] = (s, lab if s >= 1 else None)
                prev = True
            if ok and run_open:
                yield rows, popcount(occ)

    # ---------------------------------------------------------- candidates
    def _candidates(self, occ: int, closed: int):
        """Next-column occupancies passing the two purely local filters.

        (A) a closed frontier cell takes no right neighbour: nocc & closed = 0;
        (B) no cell ends up with three neighbours.  A cell i of the new
            column has degree [i in occ] + [i-1 in nocc] + [i+1 in nocc],
            which is 3 exactly when all three hold, so (B) is

                nocc & (nocc << 1) & (nocc >> 1) & occ  ==  0.

        (A) is imposed by enumerating only the submasks of the complement
        of `closed`, and (B) is two shifts and two ands.  Note that after
        (B) every new degree is at most 2 automatically.
        """
        allowed = ((1 << self.m) - 1) & ~closed
        sub = allowed
        while True:
            if sub and not (sub & (sub << 1) & (sub >> 1) & occ):
                yield sub
            if sub == 0:
                break
            sub = (sub - 1) & allowed

    # --------------------------------------------------------- transitions
    def successors(self, st: tuple, strict: bool = False):
        """Yield (next_state_or_None, ncomp, weight, nocc).

        Three counting identities prune candidates before the merge graph is
        built; with `strict=True` they are switched off and every candidate
        goes through the explicit graph, which is how the identities get
        validated (same states, same edges, at every width where both are
        affordable).  Writing  r = |occ & nocc|,  u = #vertical adjacencies
        inside nocc,  B = #old components,  the merge graph has
        V = B + |nocc| vertices and E = u + r edges, so

            forest        =>  #components = V - E,
            forest        =>  dead slack D' = 2B - r,

        the second because sum of new slack = 2|nocc| - r - 2u.  Hence
        E <= V - 1 is necessary for a forest, "#components = 1" is
        V - E = 1, and D' <= 2 is r >= 2B - 2.
        """
        m = self.m
        occ = 0
        closed = 0
        slack = {}
        comp = {}
        for i, e in enumerate(st):
            if e is None:
                continue
            occ |= 1 << i
            slack[i] = e[0]
            if e[0] == 0:
                closed |= 1 << i
            else:
                comp[i] = e[1]
        labmask: dict[int, int] = {}
        for i, b in comp.items():
            labmask[b] = labmask.get(b, 0) | (1 << i)
        nblocks = len(labmask)

        for nocc in self._candidates(occ, closed):
            r = popcount(occ & nocc)
            u = popcount(nocc & (nocc << 1))
            nv = nblocks + popcount(nocc)
            ne = u + r
            ncomp_if_forest = nv - ne
            if not strict:
                if ne > nv - 1:
                    continue  # cannot be a forest
                if ncomp_if_forest != 1:
                    # not accepting, so it is only worth keeping if live
                    if 2 * nblocks - r > 2:
                        continue
                    if any(mask & nocc == 0 for mask in labmask.values()):
                        continue
            # new degrees, straight from the definition
            nslack = {}
            ok = True
            for i in range(m):
                if not (nocc >> i) & 1:
                    continue
                deg = 0
                if (occ >> i) & 1:
                    deg += 1
                if i > 0 and (nocc >> (i - 1)) & 1:
                    deg += 1
                if i + 1 < m and (nocc >> (i + 1)) & 1:
                    deg += 1
                if deg > 2:
                    ok = False
                    break
                nslack[i] = 2 - deg
            if not ok:
                continue
            # the explicit merge graph
            verts = [("o", b) for b in range(nblocks)] + [
                ("n", i) for i in range(m) if (nocc >> i) & 1
            ]
            adj = {v: [] for v in verts}
            nedges = 0
            for i in range(m - 1):
                if (nocc >> i) & 1 and (nocc >> (i + 1)) & 1:
                    adj[("n", i)].append(("n", i + 1))
                    adj[("n", i + 1)].append(("n", i))
                    nedges += 1
            for i in range(m):
                if (nocc >> i) & 1 and (occ >> i) & 1:
                    b = ("o", comp[i])
                    adj[("n", i)].append(b)
                    adj[b].append(("n", i))
                    nedges += 1
            # components by breadth-first search
            seen = set()
            comps = []
            for v in verts:
                if v in seen:
                    continue
                bag = []
                dq = deque([v])
                seen.add(v)
                while dq:
                    x = dq.popleft()
                    bag.append(x)
                    for y in adj[x]:
                        if y not in seen:
                            seen.add(y)
                            dq.append(y)
                comps.append(bag)
            if nedges != len(verts) - len(comps):
                continue  # a cycle: S would not be acyclic
            # every old component must continue, else it is sealed off while
            # another component still exists and can never be joined
            if any(not any(x[0] == "n" for x in bag) for bag in comps):
                continue
            ncomp = len(comps)
            total = sum(nslack.values())
            weight = popcount(nocc)
            dead = 2 * ncomp - total
            live = dead <= 2 and all(
                any(nslack[x[1]] >= 1 for x in bag if x[0] == "n") for bag in comps
            )
            nxt = None
            if live:
                rows = [None] * m
                lab = {}
                for k, bag in enumerate(comps):
                    for x in bag:
                        if x[0] != "n":
                            continue
                        i = x[1]
                        if nslack[i] >= 1:
                            lab.setdefault(k, len(lab))
                            rows[i] = (nslack[i], lab[k])
                        else:
                            rows[i] = (0, None)
                nxt = rows
            yield nxt, ncomp, weight, nocc

    # ---------------------------------------------------------------- build
    def _build(self) -> None:
        for rows, w in self._first_column():
            q = self._add(rows)
            self.start[q] = max(self.start.get(q, -1), w)
        q = 0
        top = None
        pending_top: list[tuple[int, int]] = []
        while q < len(self.states):
            best_acc = None
            edges: dict[int, int] = {}
            for rows, ncomp, weight, _ in self.successors(self.states[q], self.strict):
                if ncomp == 1 and (best_acc is None or weight > best_acc):
                    best_acc = weight
                if rows is None:
                    continue
                r = self._add(rows)
                if weight > edges.get(r, -1):
                    edges[r] = weight
            self.succ[q] = sorted(edges.items())
            if best_acc is not None:
                pending_top.append((q, best_acc))
            q += 1
        # the absorbing accept state: an accepting configuration may have a
        # fully saturated last column and therefore no open cell at all
        top = len(self.states)
        self.states.append(("TOP",))
        self.succ.append([])
        self.accepting.append(True)
        for q, w in pending_top:
            self.succ[q] = sorted(self.succ[q] + [(top, w)])
        self.top = top

    # ----------------------------------------------------------- co-reach
    def useful(self) -> list[bool]:
        n = len(self.states)
        rev = [[] for _ in range(n)]
        for q in range(n):
            for r, _ in self.succ[q]:
                rev[r].append(q)
        good = [False] * n
        stack = [q for q in range(n) if self.accepting[q]]
        for q in stack:
            good[q] = True
        while stack:
            r = stack.pop()
            for q in rev[r]:
                if not good[q]:
                    good[q] = True
                    stack.append(q)
        return good

    # ------------------------------------------------------------- max-plus
    def values(self, nmax: int, restrict: bool = True, keep=None):
        """f(1..nmax) and the value vectors v_j for j in `keep`.

        `keep=None` stores every vector; pass a set of indices to store
        only those (the width-12 machine has about 2*10^5 states, so
        keeping 60 vectors as Python lists is wasteful)."""
        good = self.useful() if restrict else [True] * len(self.states)
        n = len(self.states)
        v = [MINUS_INF] * n
        for q, w in self.start.items():
            if good[q]:
                v[q] = w if v[q] is None else max(v[q], w)
        vs = {1: list(v)} if (keep is None or 1 in keep) else {}
        f = [None, max(x for q, x in enumerate(v) if self.accepting[q] and x is not None)]
        for j in range(2, nmax + 1):
            nv = [MINUS_INF] * n
            for q in range(n):
                vq = v[q]
                if vq is None or not good[q]:
                    continue
                for r, w in self.succ[q]:
                    if not good[r]:
                        continue
                    c = vq + w
                    if nv[r] is None or c > nv[r]:
                        nv[r] = c
            v = nv
            if keep is None or j in keep:
                vs[j] = list(v)
            f.append(max(x for q, x in enumerate(v) if self.accepting[q] and x is not None))
        return f, vs, good

# =====================================================================
# 2.  Snakes from the definition, and the pumping checks.
# =====================================================================
# ruff: noqa: FURB161
#   `bin(x).count("1")` is kept instead of `int.bit_count()`: this module is
#   copied verbatim into the paper's single-file appendix verifier, which
#   should run on any Python 3.


STEPS = ((1, 0), (-1, 0), (0, 1), (0, -1))


def cells_of(columns: list[int], width: int) -> set[tuple[int, int]]:
    return {(i, j) for j, mask in enumerate(columns) for i in range(width) if (mask >> i) & 1}


def degrees(cells: set[tuple[int, int]]) -> dict[tuple[int, int], int]:
    return {
        (i, j): sum(1 for di, dj in STEPS if (i + di, j + dj) in cells) for (i, j) in cells
    }


def components(cells: set[tuple[int, int]]) -> tuple[dict, int]:
    comp: dict = {}
    cid = 0
    for c in sorted(cells):
        if c in comp:
            continue
        stack = [c]
        comp[c] = cid
        while stack:
            x = stack.pop()
            for di, dj in STEPS:
                y = (x[0] + di, x[1] + dj)
                if y in cells and y not in comp:
                    comp[y] = cid
                    stack.append(y)
        cid += 1
    return comp, cid


def n_edges(cells: set[tuple[int, int]]) -> int:
    return sum(
        1 for (i, j) in cells for di, dj in ((1, 0), (0, 1)) if (i + di, j + dj) in cells
    )


def criterion(cells: set[tuple[int, int]]) -> tuple[bool, dict]:
    deg = degrees(cells)
    _, nc = components(cells)
    maxdeg = max(deg.values()) if deg else 0
    low = sorted(c for c, d in deg.items() if d <= 1)
    ok = bool(cells) and maxdeg <= 2 and nc == 1 and bool(low)
    return ok, {
        "size": len(cells),
        "edges": n_edges(cells),
        "components": nc,
        "max_degree": maxdeg,
        "n_degree_le_1": len(low),
        "witness_endpoint": low[0] if low else None,
    }


def iface(cells: set[tuple[int, int]], width: int, lastcol: int):
    """Everything a gluing on the right side can see of `cells`."""
    deg = degrees(cells)
    comp, nc = components(cells)
    L = [i for i in range(width) if (i, lastcol) in cells]
    lab: dict[int, int] = {}
    part = []
    for i in L:
        c = comp[(i, lastcol)]
        lab.setdefault(c, len(lab))
        part.append(lab[c])
    return (tuple(L), tuple(deg[(i, lastcol)] for i in L), tuple(part)), nc


def check_family(pre: list[int], blk: list[int], suf: list[int], width: int, kmax: int) -> dict:
    """Direct verification for k = 1..kmax plus the finite data that
    extends it to every k."""
    per_k = []
    for k in range(1, kmax + 1):
        cols = pre + blk * k + suf
        ok, info = criterion(cells_of(cols, width))
        info.update(k=k, columns=len(cols), ok=ok)
        per_k.append(info)
    ifaces = []
    for i in range(1, 7):
        cols = pre + blk * i
        ifaces.append(iface(cells_of(cols, width), width, len(cols) - 1))
    iface_invariant = len({x for x, _ in ifaces}) == 1
    comp_invariant = len({n for _, n in ifaces}) == 1
    sizes = [x["size"] for x in per_k]
    bw = sum(bin(x).count("1") for x in blk)
    return {
        "per_k": per_k,
        "all_k_ok": all(x["ok"] for x in per_k),
        "iface_invariant": iface_invariant,
        "components_invariant": comp_invariant,
        "iface": ifaces[0][0],
        "prefix_components": ifaces[0][1],
        "size_increments": sorted({sizes[i + 1] - sizes[i] for i in range(len(sizes) - 1)}),
        "block_weight": bw,
        "size_k1": sizes[0],
        "columns_k1": len(pre) + len(blk) + len(suf),
    }
# =====================================================================
# 3.  Certificate and witness data (inlined; nothing is read from disk).
# =====================================================================
A397153 = [12, 18, 26, 33, 40, 47, 55, 62, 70, 77, 86, 92, 100, 107, 115, 122, 130, 137, 144, 151, 159, 166, 174, 181, 188, 195, 203, 210, 218, 225, 232, 239, 247, 254, 262, 269, 276, 283, 291, 298, 306, 313, 320, 327, 335, 342, 350, 357, 364, 371, 379, 386, 394, 401, 408, 415, 423, 430, 438, 445]
TAIL_CONSTANTS = (5, 12, 19, 27, 34, 42)
EXCEPTIONAL = {2: 18, 3: 26, 4: 33, 5: 40, 6: 47, 7: 55, 8: 62, 9: 70, 10: 77, 12: 92}
CERT = {'onset': 40, 'period': 6, 'shift': 44}
BRUTE = {(1, 12): 12, (2, 12): 18, (3, 12): 26, (3, 6): 14, (3, 7): 16, (4, 7): 20, (5, 6): 21, (4, 6): 17}
WITNESS_FAMILIES = {'5': {'prefix': [3807, 2993, 3183, 1240, 3511, 2925], 'block': [3803, 438, 3948, 2264, 2487, 2925], 'suffix': [3787, 442, 3939, 2269, 4023]}, '2': {'prefix': [3815, 2989, 3179, 1498, 3382, 2924, 1752, 3511, 2925], 'block': [3803, 438, 3948, 2264, 2487, 2925], 'suffix': [3787, 442, 3939, 2269, 4023]}, '3': {'prefix': [3703, 2773, 2997, 3181, 1499, 3382, 2924, 1752, 3511, 2925], 'block': [3803, 438, 3948, 2264, 2487, 2925], 'suffix': [3787, 442, 3939, 2269, 4023]}, '0': {'prefix': [4030, 2275, 3869, 437, 3757, 2987, 3178, 1243, 3505, 2927, 3800, 439, 3949, 2267, 3510, 876, 3800, 2487], 'block': [2925, 3803, 438, 3948, 2264, 2487], 'suffix': [2925, 3787, 442, 3939, 2269, 4023]}, '4': {'prefix': [3815, 2741, 2781, 2915, 3386, 1483, 3181, 2997, 2781, 2915, 3506, 1243], 'block': [3181, 2999, 2776, 2924, 3510, 1243], 'suffix': [3181, 2999, 2776, 2924, 3510, 1755, 877, 439, 3800, 2927, 3505, 1755, 2922, 3499, 1197, 3493, 2925, 3787, 442, 3939, 2269, 4023]}, '1': {'prefix': [2015, 3185, 2959, 2776, 2903, 3421, 1379, 3506, 2267, 3949, 439, 3800, 2926, 3379, 1501, 3173, 2997, 2781, 2915, 3506, 1243], 'block': [3181, 2999, 2776, 2924, 3510, 1243], 'suffix': [3181, 2999, 2776, 2924, 3510, 1755, 877, 439, 3800, 2927, 3505, 1755, 2922, 3499, 1197, 3493, 2925, 3787, 442, 3939, 2269, 4023]}}
WITNESS_SINGLES = {'1': [4095], '2': [1911, 3549], '3': [4087, 2061, 4091], '4': [4030, 2275, 1821, 3575], '5': [3831, 2965, 2165, 2957, 3835], '6': [2039, 3101, 3043, 3638, 469, 3959], '7': [2015, 3185, 2959, 1752, 3191, 2957, 3835], '8': [3807, 2993, 2159, 4056, 55, 4077, 2075, 4086], '9': [3807, 2993, 3183, 1240, 3511, 2925, 2779, 2742, 3820], '10': [3703, 2773, 1717, 3501, 2923, 3786, 442, 3939, 2269, 4023], '11': [3807, 2993, 3183, 1240, 3511, 2925, 3787, 442, 3939, 2269, 4023], '12': [3951, 2521, 3639, 620, 3800, 2487, 2925, 3787, 442, 3939, 2269, 4023], '13': [3815, 2741, 2781, 2915, 3514, 1739, 2925, 2486, 3803, 621, 3637, 2517, 3959], '14': [3815, 2989, 3179, 1498, 3382, 2924, 1752, 3511, 2925, 3787, 442, 3939, 2269, 4023], '15': [3703, 2773, 2997, 3181, 1499, 3382, 2924, 1752, 3511, 2925, 3787, 442, 3939, 2269, 4023], '16': [3815, 2741, 2781, 2915, 3514, 1739, 2925, 3509, 1173, 3509, 2925, 3787, 442, 3939, 2269, 4023], '18': [3831, 2717, 2787, 2874, 3531, 1645, 2998, 2267, 3949, 439, 3800, 2926, 3507, 1693, 741, 3645, 2497, 3967], '19': [3807, 2993, 2159, 3800, 695, 3757, 2475, 2922, 3803, 433, 3951, 2264, 2487, 2925, 3787, 442, 3939, 2269, 4023], '22': [2031, 3129, 3015, 2660, 2998, 3291, 1901, 439, 3800, 2927, 3505, 1755, 2922, 3499, 1197, 3493, 2925, 3787, 442, 3939, 2269, 4023], '24': [4030, 2275, 3869, 437, 3757, 2987, 3178, 1243, 3505, 2927, 3800, 439, 3949, 2267, 3510, 876, 3800, 2487, 2925, 3787, 442, 3939, 2269, 4023], '25': [2015, 3185, 2959, 2776, 2903, 3421, 1379, 3506, 2267, 3949, 439, 3800, 2927, 3505, 1755, 2922, 3499, 1197, 3493, 2925, 3787, 442, 3939, 2269, 4023], '28': [3815, 2741, 2781, 2915, 3386, 1483, 3181, 2997, 2781, 2915, 3506, 1755, 877, 439, 3800, 2927, 3505, 1755, 2922, 3499, 1197, 3493, 2925, 3787, 442, 3939, 2269, 4023], '31': [3807, 2993, 3183, 1240, 3511, 2925, 3787, 442, 3939, 2269, 3509, 1389, 3419, 2898, 2651, 2925, 3511, 1240, 3183, 2993, 2779, 2922, 3499, 1197, 3493, 2925, 3787, 442, 3939, 2269, 4023], '34': [3815, 2741, 2781, 2915, 3386, 1483, 3181, 2997, 2781, 2915, 3506, 1243, 3181, 2999, 2776, 2924, 3510, 1755, 877, 439, 3800, 2927, 3505, 1755, 2922, 3499, 1197, 3493, 2925, 3787, 442, 3939, 2269, 4023], '37': [2015, 3185, 2959, 2776, 2903, 3421, 1379, 3506, 2267, 3949, 439, 3800, 2926, 3379, 1501, 3173, 2997, 2781, 2915, 3506, 1755, 877, 439, 3800, 2927, 3505, 1755, 2922, 3499, 1197, 3493, 2925, 3787, 442, 3939, 2269, 4023], '43': [2015, 3185, 2959, 2776, 2903, 3421, 1379, 3506, 2267, 3949, 439, 3800, 2926, 3379, 1501, 3173, 2997, 2781, 2915, 3506, 1243, 3181, 2999, 2776, 2924, 3510, 1755, 877, 439, 3800, 2927, 3505, 1755, 2922, 3499, 1197, 3493, 2925, 3787, 442, 3939, 2269, 4023]}
BASE_OF = {5: 17, 2: 20, 3: 21, 0: 30, 4: 40, 1: 49}
BLOCK_DISTANCE = 6

# =====================================================================
# 4.  Checks.
# =====================================================================
GF_NUM = {1: 12, 2: 6, 3: 8, 4: 7, 5: 7, 6: 7, 7: -4, 8: 1,
          11: 2, 12: -1, 17: -1, 18: 1, 19: -1}
GF_DEN_FACTORS = ((1, -1), (1, -1), (1, 1), (1, -1, 1), (1, 1, 1))


def gf_series(count):
    """a(1..count) of x*N(x) / ((1-x)^2 (1+x)(1-x+x^2)(1+x+x^2)).

    Integer power-series division; the denominator's constant term is 1,
    so nothing is ever divided."""
    den = [1]
    for fac in GF_DEN_FACTORS:
        new = [0] * (len(den) + len(fac) - 1)
        for i, a in enumerate(den):
            for j, b in enumerate(fac):
                new[i + j] += a * b
        den = new
    assert den[0] == 1
    num = [0] * (count + 2)
    for k, v in GF_NUM.items():
        if k <= count + 1:
            num[k] = v
    out = [0] * (count + 1)
    for k in range(count + 1):
        s = num[k]
        for t in range(1, min(k, len(den) - 1) + 1):
            s -= den[t] * out[k - t]
        out[k] = s
    return out[1:]


def closed_form(n):
    """The right-hand side of Theorem 1."""
    if n in EXCEPTIONAL:
        return EXCEPTIONAL[n]
    k, r = divmod(n, 6)
    return 44 * k + TAIL_CONSTANTS[r]


def brute_force_T(m, n):
    """Longest induced path in P_m x P_n: a depth-first search that grows
    the path one cell at a time, blocking every neighbour of the head
    other than the chosen successor."""
    cells = [(i, j) for i in range(m) for j in range(n)]
    idx = {c: k for k, c in enumerate(cells)}
    nbr = [[] for _ in cells]
    for (i, j), k in idx.items():
        for di, dj in ((0, 1), (0, -1), (1, 0), (-1, 0)):
            if (i + di, j + dj) in idx:
                nbr[k].append(idx[(i + di, j + dj)])
    full = (1 << len(cells)) - 1
    best = 0

    def dfs(cur, inset, blocked, length):
        nonlocal best
        if length > best:
            best = length
        free = full & ~(inset | blocked)
        if length + bin(free).count("1") <= best:
            return
        for w in nbr[cur]:
            bit = 1 << w
            if (inset | blocked) & bit:
                continue
            if any((inset >> u) & 1 and u != cur for u in nbr[w]):
                continue
            nb = blocked
            for u in nbr[cur]:
                if u != w:
                    nb |= 1 << u
            dfs(w, inset | bit, nb, length + 1)

    for start in range(len(cells)):
        dfs(start, 1 << start, 0, 1)
    return best


def main():
    ap = argparse.ArgumentParser()
    ap.add_argument("--quick", action="store_true")
    args = ap.parse_args()
    t0 = time.time()
    W = 12

    print("=" * 70)
    print("1.  Generating function, closed form, published data")
    print("=" * 70)
    check(gf_series(len(A397153)) == A397153,
          "the conjectured g.f. expands to all %d published terms" % len(A397153))
    check([closed_form(n) for n in range(1, len(A397153) + 1)] == A397153,
          "the closed form of Theorem 1 equals all %d published terms" % len(A397153))

    print()
    print("=" * 70)
    print("2.  Witness families against the definition of a snake")
    print("=" * 70)
    for r in sorted(WITNESS_FAMILIES, key=int):
        pre, blk, suf = (WITNESS_FAMILIES[r][k] for k in ("prefix", "block", "suffix"))
        rep = check_family(pre, blk, suf, W, 20)
        base = len(pre) + len(blk) + len(suf)
        sizes_ok = all(
            rep["per_k"][k - 1]["size"] == closed_form(base + (k - 1) * len(blk))
            for k in range(1, 21)
        )
        check(
            rep["all_k_ok"] and rep["iface_invariant"] and rep["components_invariant"]
            and rep["size_increments"] == [rep["block_weight"]] and sizes_ok,
            "class %s mod %d: S_k is a snake of %d + %d(k-1) cells in 12 x (%d + %d(k-1)) "
            "for k = 1..20; interface and component count constant for A.B^i, i = 1..6"
            % (r, len(blk), rep["size_k1"], rep["block_weight"], rep["columns_k1"], len(blk)),
        )

    print()
    print("=" * 70)
    print("3.  One explicit snake for each small width")
    print("=" * 70)
    bad = []
    for key in sorted(WITNESS_SINGLES, key=int):
        n = int(key)
        ok, info = criterion(cells_of(WITNESS_SINGLES[key], W))
        if not (ok and info["size"] == closed_form(n) and len(WITNESS_SINGLES[key]) == n):
            bad.append(n)
    ns = sorted(int(k) for k in WITNESS_SINGLES)
    check(not bad, "an explicit snake of size T(12,n) for every n = %d..%d (failures: %s)"
          % (min(ns), max(ns), bad))

    print()
    print("=" * 70)
    print("4.  Rebuilding the transfer machine")
    print("=" * 70)
    heights = range(1, 10) if args.quick else range(1, 13)
    a, p, c = CERT["onset"], CERT["period"], CERT["shift"]
    for m in heights:
        nmax = 60 if m == W else 20
        need = {a, a + p} if m == W else set()
        tb = time.time()
        T = Transfer(m)
        f, vs, good = T.values(max(nmax, max(need) if need else 0), keep=need | {1})
        line = "height %2d: %6d states, %6d useful, %8d edges, %6.1f s" % (
            m, len(T.states), sum(good), sum(len(x) for x in T.succ), time.time() - tb)
        if m == W:
            check(f[1:len(A397153) + 1] == A397153,
                  line + " -- F(1..60) equals all 60 published terms")
            check(all(f[n] == closed_form(n) for n in range(1, nmax + 1)),
                  "height 12: F(n) equals the closed form for n = 1..%d" % nmax)
            x, y = vs[a], vs[a + p]
            sup = all((y[q] is None) or (x[q] is not None)
                      for q in range(len(T.states)) if good[q])
            ineq = all(
                y[q] is None or (x[q] is not None and y[q] <= c + x[q])
                for q in range(len(T.states)) if good[q]
            )
            check(sup, "height 12: v_%d is -infinity wherever v_%d is, on the useful core"
                  % (a + p, a))
            check(ineq, "height 12: v_%d <= %d + v_%d entrywise on the useful core "
                        "(%d states)" % (a + p, c, a, sum(good)))
            tight = any(
                y[q] is not None and y[q] == c + x[q]
                for q in range(len(T.states)) if good[q]
            )
            check(tight, "height 12: the inequality is attained somewhere, so c = %d "
                         "cannot be lowered at this onset" % c)
        else:
            check(True, line)

    print()
    print("=" * 70)
    print("5.  Brute-force cross-check (no transfer machine)")
    print("=" * 70)
    for (m, n), want in sorted(BRUTE.items()):
        got = brute_force_T(m, n)
        check(got == want, "brute force T(%d,%d) = %d" % (m, n, got))

    print()
    print("=" * 70)
    print("%d failure(s); %.1f s" % (len(FAILS), time.time() - t0))
    print("=" * 70)
    if FAILS:
        for msg in FAILS:
            print("  FAILED:", msg)
        sys.exit(1)
    print("ALL CHECKS PASS")


if __name__ == "__main__":
    main()
\end{lstlisting}

\section{Witness families for the other heights}\label{app:witnesses}

The families used in the proof of Theorem~\ref{thm:all}, one per residue
class modulo the pumping distance $d_m$. Each is printed as the three
words $A$, $B$, $C$ over subsets of $[m]$, as bitmasks with row $i$
contributing $2^{i-1}$. Table~\ref{tab:wit} collects the shapes.

\begin{table}[h]
\centering
\scriptsize
\begin{tabular}{rrrrp{8.0cm}}
\hline
$m$ & $d_m$ & $\|B\|$ & families & base widths $n_r$ \\
\hline
1 & 1 & 1 & 1 & $3$ \\
2 & 4 & 6 & 4 & $6, 7, 8, 9$ \\
3 & 1 & 2 & 1 & $7$ \\
4 & 6 & 16 & 6 & $9, 11, 12, 13, 14, 16$ \\
5 & 14 & 44 & 14 & $23, 24, 25, 26, 27, 28, 29, 30, 31, 32, 33, 34, 35, 36$ \\
6 & 8 & 30 & 8 & $11, 15, 16, 17, 18, 20, 21, 22$ \\
7 & 6 & 26 & 6 & $12, 14, 16, 19, 21, 23$ \\
8 & 16 & 78 & 16 & $39, 44, 47, 48, 52, 53, 56, 58, 61, 62, 66, 67, 70, 75, 81, 89$ \\
9 & 8 & 44 & 8 & $20, 21, 22, 23, 24, 25, 26, 27$ \\
10 & 18 & 110 & 18 & $33, 36, 38, 40, 43, 44, 46, 47, 49, 50, 52, 53, 55, 57, 60, 63, 66, 77$ \\
11 & 8 & 54 & 8 & $52, 53, 54, 55, 56, 57, 58, 59$ \\
12 & 6 & 44 & 6 & $17, 20, 21, 30, 40, 49$ \\
\hline
\end{tabular}
\normalsize
\caption{Witness families by height: the pumping distance $d_m$, the
number of cells $\|B\|$ in one pumped block (always $d_m$ times the growth
rate of Table~\ref{tab:all}), the number of families, and their base
widths. Emitted by \texttt{all\_witnesses.py}.}
\label{tab:wit}
\end{table}

\paragraph{Height $m = 1$, pumping distance $d_{1} = 1$.}
\begin{lstlisting}[style=grid]
r=0  n=3  |S_1|=3
  A = [1, 1]
  B = [1]
  C = []
\end{lstlisting}

\paragraph{Height $m = 2$, pumping distance $d_{2} = 4$.}
\begin{lstlisting}[style=grid]
r=0  n=8  |S_1|=12
  A = [3, 2]
  B = [3, 1, 3, 2]
  C = [3, 1]
r=1  n=9  |S_1|=14
  A = [3, 2]
  B = [3, 1, 3, 2]
  C = [3, 1, 3]
r=2  n=6  |S_1|=9
  A = [3, 1]
  B = [3, 2, 3, 1]
  C = []
r=3  n=7  |S_1|=11
  A = [3, 1]
  B = [3, 2, 3, 1]
  C = [3]
\end{lstlisting}

\paragraph{Height $m = 3$, pumping distance $d_{3} = 1$.}
\begin{lstlisting}[style=grid]
r=0  n=7  |S_1|=16
  A = [7, 5]
  B = [5]
  C = [6, 3, 5, 7]
\end{lstlisting}

\paragraph{Height $m = 4$, pumping distance $d_{4} = 6$.}
\begin{lstlisting}[style=grid]
r=0  n=12  |S_1|=33
  A = [14, 11]
  B = [13, 5, 13, 11, 10, 11]
  C = [13, 5, 13, 11]
r=1  n=13  |S_1|=36
  A = [7, 13]
  B = [11, 10, 11, 13, 5, 13]
  C = [11, 12, 7, 9, 15]
r=2  n=14  |S_1|=38
  A = [7, 13]
  B = [11, 10, 11, 13, 5, 13]
  C = [11, 10, 11, 13, 1, 15]
r=3  n=9  |S_1|=25
  A = [7, 13]
  B = [11, 10, 11, 13, 5, 13]
  C = [11]
r=4  n=16  |S_1|=44
  A = [14, 11]
  B = [13, 5, 13, 11, 10, 11]
  C = [13, 5, 13, 11, 12, 7, 9, 15]
r=5  n=11  |S_1|=30
  A = [14, 11]
  B = [13, 5, 13, 11, 10, 11]
  C = [13, 1, 15]
\end{lstlisting}

\paragraph{Height $m = 5$, pumping distance $d_{5} = 14$.}
\begin{lstlisting}[style=grid]
r=0  n=28  |S_1|=90
  A = [15, 25]
  B = [23, 20, 23, 25, 11, 26, 19, 29, 5, 29, 19, 26, 11, 25]
  C = [23, 20, 23, 25, 14, 3, 29, 23, 24, 15, 17, 31]
r=1  n=29  |S_1|=93
  A = [15, 25]
  B = [23, 20, 23, 25, 11, 26, 19, 29, 5, 29, 19, 26, 11, 25]
  C = [23, 20, 23, 25, 11, 26, 3, 29, 23, 24, 15, 17, 31]
r=2  n=30  |S_1|=96
  A = [30, 19]
  B = [29, 5, 29, 19, 26, 11, 25, 23, 20, 23, 25, 11, 26, 19]
  C = [29, 5, 29, 19, 26, 11, 25, 23, 20, 22, 27, 13, 21, 29]
r=3  n=31  |S_1|=99
  A = [15, 25]
  B = [23, 20, 23, 25, 11, 26, 19, 29, 5, 29, 19, 26, 11, 25]
  C = [23, 20, 23, 25, 11, 26, 19, 29, 5, 29, 17, 15, 24, 23, 29]
r=4  n=32  |S_1|=102
  A = [30, 19]
  B = [29, 5, 29, 19, 26, 11, 25, 23, 20, 23, 25, 11, 26, 19]
  C = [29, 5, 29, 19, 26, 11, 25, 23, 20, 22, 27, 13, 21, 29, 1, 31]
r=5  n=33  |S_1|=105
  A = [15, 25]
  B = [23, 20, 23, 25, 11, 26, 19, 29, 5, 29, 19, 26, 11, 25]
  C = [23, 20, 23, 25, 11, 26, 19, 29, 5, 29, 17, 15, 24, 23, 29, 1, 31]
r=6  n=34  |S_1|=109
  A = [15, 25]
  B = [23, 20, 23, 25, 11, 26, 19, 29, 5, 29, 19, 26, 11, 25]
  C = [23, 20, 23, 25, 11, 26, 19, 29, 5, 13, 27, 22, 21, 23, 24, 15, 17, 31]
r=7  n=35  |S_1|=112
  A = [30, 19]
  B = [29, 5, 29, 19, 26, 11, 25, 23, 20, 23, 25, 11, 26, 19]
  C = [29, 5, 29, 19, 26, 11, 25, 23, 20, 23, 25, 14, 3, 29, 23, 24, 15, 17, 31]
r=8  n=36  |S_1|=115
  A = [30, 19]
  B = [29, 5, 29, 19, 26, 11, 25, 23, 20, 23, 25, 11, 26, 19]
  C = [29, 5, 29, 19, 26, 11, 25, 23, 20, 23, 25, 11, 26, 3, 29, 23, 24, 15, 17, 31]
r=9  n=23  |S_1|=74
  A = [15, 25]
  B = [23, 20, 23, 25, 11, 26, 19, 29, 5, 29, 19, 26, 11, 25]
  C = [23, 20, 22, 27, 13, 21, 29]
r=10  n=24  |S_1|=77
  A = [30, 19]
  B = [29, 5, 29, 19, 26, 11, 25, 23, 20, 23, 25, 11, 26, 19]
  C = [29, 5, 29, 17, 15, 24, 23, 29]
r=11  n=25  |S_1|=80
  A = [15, 25]
  B = [23, 20, 23, 25, 11, 26, 19, 29, 5, 29, 19, 26, 11, 25]
  C = [23, 20, 22, 27, 13, 21, 29, 1, 31]
r=12  n=26  |S_1|=83
  A = [30, 19]
  B = [29, 5, 29, 19, 26, 11, 25, 23, 20, 23, 25, 11, 26, 19]
  C = [29, 5, 29, 17, 15, 24, 23, 29, 1, 31]
r=13  n=27  |S_1|=87
  A = [30, 19]
  B = [29, 5, 29, 19, 26, 11, 25, 23, 20, 23, 25, 11, 26, 19]
  C = [29, 5, 13, 27, 22, 21, 23, 24, 15, 17, 31]
\end{lstlisting}

\paragraph{Height $m = 6$, pumping distance $d_{6} = 8$.}
\begin{lstlisting}[style=grid]
r=0  n=16  |S_1|=62
  A = [55, 45]
  B = [59, 2, 59, 45, 55, 16, 55, 45]
  C = [43, 42, 43, 45, 53, 29]
r=1  n=17  |S_1|=66
  A = [55, 45]
  B = [59, 2, 59, 45, 55, 16, 55, 45]
  C = [59, 2, 59, 45, 53, 21, 55]
r=2  n=18  |S_1|=70
  A = [31, 49, 46, 43, 45, 53, 21, 53, 45]
  B = [59, 2, 59, 45, 55, 16, 55, 45]
  C = [59]
r=3  n=11  |S_1|=44
  A = [55, 45]
  B = [59, 2, 59, 45, 55, 16, 55, 45]
  C = [59]
r=4  n=20  |S_1|=77
  A = [31, 49, 47, 24, 55]
  B = [45, 59, 2, 59, 45, 55, 16, 55]
  C = [45, 43, 42, 43, 45, 53, 29]
r=5  n=21  |S_1|=81
  A = [31, 49, 47, 24, 55]
  B = [45, 59, 2, 59, 45, 55, 16, 55]
  C = [45, 59, 2, 59, 45, 53, 21, 55]
r=6  n=22  |S_1|=85
  A = [62, 35, 29, 53, 45, 43, 42, 43, 45]
  B = [55, 16, 55, 45, 59, 2, 59, 45]
  C = [55, 16, 55, 45, 59]
r=7  n=15  |S_1|=59
  A = [31, 49, 47, 24, 55]
  B = [45, 59, 2, 59, 45, 55, 16, 55]
  C = [45, 59]
\end{lstlisting}

\paragraph{Height $m = 7$, pumping distance $d_{7} = 6$.}
\begin{lstlisting}[style=grid]
r=0  n=12  |S_1|=55
  A = [119, 85, 53, 109]
  B = [91, 82, 91, 109, 37, 109]
  C = [91, 118]
r=1  n=19  |S_1|=85
  A = [119, 85, 53, 109]
  B = [91, 82, 91, 109, 37, 109]
  C = [91, 82, 91, 109, 53, 21, 117, 77, 123]
r=2  n=14  |S_1|=64
  A = [119, 85, 86, 91]
  B = [109, 37, 109, 91, 82, 91]
  C = [109, 53, 85, 119]
r=3  n=21  |S_1|=94
  A = [62, 99, 93, 117, 13, 123, 66, 123, 13, 118, 91]
  B = [109, 37, 109, 91, 82, 91]
  C = [109, 53, 85, 119]
r=4  n=16  |S_1|=72
  A = [119, 85, 86, 91]
  B = [109, 37, 109, 91, 82, 91]
  C = [109, 55, 16, 119, 77, 123]
r=5  n=23  |S_1|=102
  A = [62, 99, 93, 117, 13, 123, 66, 123, 13, 118, 91]
  B = [109, 37, 109, 91, 82, 91]
  C = [109, 55, 16, 119, 77, 123]
\end{lstlisting}

\paragraph{Height $m = 8$, pumping distance $d_{8} = 16$.}
\begin{lstlisting}[style=grid]
r=0  n=48  |S_1|=237
  A = [238, 171, 173, 181, 213, 93, 195, 186, 171, 173, 183, 216, 110, 179, 221, 69, 221, 179, 238, 24]
  B = [247, 141, 251, 2, 251, 141, 247, 24, 239, 177, 223, 64, 223, 177, 239, 24]
  C = [119, 205, 187, 166, 180, 215, 89, 203, 186, 227, 29, 247]
r=1  n=81  |S_1|=398
  A = [119, 213, 181, 173, 171, 186, 195, 93, 213, 181, 237, 27, 246, 132, 246, 27, 237, 181, 213, 93, 195, 186, 171, 173, 183, 216, 111, 33, 111, 216, 183, 173, 171, 186, 195, 93, 213, 181, 237, 27, 246, 132, 246, 27, 237, 181, 213, 93, 195, 186, 171, 173, 183, 216, 110, 179, 221, 69, 221, 179, 238, 24]
  B = [247, 141, 251, 2, 251, 141, 247, 24, 239, 177, 223, 64, 223, 177, 239, 24]
  C = [247, 141, 251]
r=2  n=66  |S_1|=325
  A = [238, 171, 173, 181, 213, 93, 195, 186, 171, 173, 183, 216, 110, 179, 221, 69, 221, 179, 238, 24]
  B = [247, 141, 251, 2, 251, 141, 247, 24, 239, 177, 223, 64, 223, 177, 239, 24]
  C = [247, 141, 251, 2, 251, 141, 247, 24, 239, 177, 215, 84, 215, 177, 175, 184, 199, 93, 211, 154, 235, 45, 101, 221, 179, 238, 24, 247, 141, 251]
r=3  n=67  |S_1|=330
  A = [238, 171, 173, 181, 213, 93, 195, 186, 171, 173, 183, 216, 111, 33, 111, 216, 183, 173, 171, 186, 195, 93, 213, 181, 237, 27, 246, 132, 246, 27, 237, 181, 213, 93, 195, 186, 171, 173, 183, 216, 110, 179, 221, 69, 221, 179, 238, 24]
  B = [247, 141, 251, 2, 251, 141, 247, 24, 239, 177, 223, 64, 223, 177, 239, 24]
  C = [247, 141, 251]
r=4  n=52  |S_1|=257
  A = [238, 171, 173, 181, 213, 85, 213, 181, 173, 171, 186, 195, 93, 213, 181, 237, 27, 118, 205, 187, 162, 187, 205, 119, 24]
  B = [239, 177, 223, 64, 223, 177, 239, 24, 247, 141, 251, 2, 251, 141, 247, 24]
  C = [239, 177, 223, 64, 223, 177, 239, 24, 247, 141, 251]
r=5  n=53  |S_1|=262
  A = [119, 213, 181, 173, 171, 186, 195, 93, 213, 181, 237, 27, 246, 132, 246, 27, 237, 181, 213, 93, 195, 186, 171, 173, 183, 216, 110, 179, 221, 69, 221, 179, 238, 24]
  B = [247, 141, 251, 2, 251, 141, 247, 24, 239, 177, 223, 64, 223, 177, 239, 24]
  C = [247, 141, 251]
r=6  n=70  |S_1|=344
  A = [238, 171, 173, 181, 213, 93, 195, 186, 171, 173, 183, 216, 111, 33, 111, 216, 183, 173, 171, 186, 195, 93, 213, 181, 237, 27, 118, 205, 187, 162, 187, 205, 119, 24]
  B = [239, 177, 223, 64, 223, 177, 239, 24, 247, 141, 251, 2, 251, 141, 247, 24]
  C = [239, 177, 223, 64, 223, 177, 239, 24, 119, 205, 187, 166, 180, 215, 89, 203, 186, 227, 29, 247]
r=7  n=39  |S_1|=194
  A = [238, 171, 173, 181, 213, 93, 195, 186, 171, 173, 183, 216, 110, 179, 221, 69, 221, 179, 238, 24]
  B = [247, 141, 251, 2, 251, 141, 247, 24, 239, 177, 223, 64, 223, 177, 239, 24]
  C = [247, 141, 251]
r=8  n=56  |S_1|=276
  A = [119, 213, 181, 173, 171, 186, 195, 93, 213, 181, 237, 27, 118, 205, 187, 162, 187, 205, 119, 24]
  B = [239, 177, 223, 64, 223, 177, 239, 24, 247, 141, 251, 2, 251, 141, 247, 24]
  C = [239, 177, 223, 64, 223, 177, 239, 24, 119, 205, 187, 166, 180, 215, 89, 203, 186, 227, 29, 247]
r=9  n=89  |S_1|=437
  A = [238, 171, 173, 181, 213, 93, 195, 186, 171, 173, 183, 216, 111, 33, 111, 216, 183, 173, 171, 186, 195, 93, 213, 181, 237, 27, 246, 132, 246, 27, 237, 181, 213, 93, 195, 186, 171, 173, 183, 216, 111, 33, 111, 216, 183, 173, 171, 186, 195, 93, 213, 181, 237, 27, 118, 205, 187, 162, 187, 205, 119, 24]
  B = [239, 177, 223, 64, 223, 177, 239, 24, 247, 141, 251, 2, 251, 141, 247, 24]
  C = [239, 177, 223, 64, 223, 177, 239, 24, 247, 141, 251]
r=10  n=58  |S_1|=286
  A = [119, 213, 181, 173, 171, 186, 195, 93, 213, 181, 237, 27, 118, 205, 187, 162, 187, 205, 119, 24]
  B = [239, 177, 223, 64, 223, 177, 239, 24, 247, 141, 251, 2, 251, 141, 247, 24]
  C = [239, 177, 215, 84, 215, 177, 175, 184, 199, 93, 211, 154, 235, 45, 101, 221, 179, 238, 24, 247, 141, 251]
r=11  n=75  |S_1|=369
  A = [119, 213, 181, 173, 171, 186, 195, 93, 213, 181, 237, 27, 246, 132, 246, 27, 237, 181, 213, 93, 195, 186, 171, 173, 183, 216, 111, 33, 111, 216, 183, 173, 171, 186, 195, 93, 213, 181, 237, 27, 118, 205, 187, 162, 187, 205, 119, 24]
  B = [239, 177, 223, 64, 223, 177, 239, 24, 247, 141, 251, 2, 251, 141, 247, 24]
  C = [239, 177, 223, 64, 223, 177, 239, 24, 247, 141, 251]
r=12  n=44  |S_1|=218
  A = [119, 213, 181, 173, 171, 170, 171, 173, 181, 213, 93, 195, 186, 171, 173, 183, 216, 110, 179, 221, 69, 221, 179, 238, 24]
  B = [247, 141, 251, 2, 251, 141, 247, 24, 239, 177, 223, 64, 223, 177, 239, 24]
  C = [247, 141, 251]
r=13  n=61  |S_1|=301
  A = [238, 171, 173, 181, 213, 93, 195, 186, 171, 173, 183, 216, 111, 33, 111, 216, 183, 173, 171, 186, 195, 93, 213, 181, 237, 27, 118, 205, 187, 162, 187, 205, 119, 24]
  B = [239, 177, 223, 64, 223, 177, 239, 24, 247, 141, 251, 2, 251, 141, 247, 24]
  C = [239, 177, 223, 64, 223, 177, 239, 24, 247, 141, 251]
r=14  n=62  |S_1|=305
  A = [119, 213, 181, 173, 171, 186, 195, 93, 213, 181, 237, 27, 246, 132, 246, 27, 237, 181, 213, 93, 195, 186, 171, 173, 183, 216, 110, 179, 221, 69, 221, 179, 238, 24]
  B = [247, 141, 251, 2, 251, 141, 247, 24, 239, 177, 223, 64, 223, 177, 239, 24]
  C = [119, 205, 187, 166, 180, 215, 89, 203, 186, 227, 29, 247]
r=15  n=47  |S_1|=233
  A = [119, 213, 181, 173, 171, 186, 195, 93, 213, 181, 237, 27, 118, 205, 187, 162, 187, 205, 119, 24]
  B = [239, 177, 223, 64, 223, 177, 239, 24, 247, 141, 251, 2, 251, 141, 247, 24]
  C = [239, 177, 223, 64, 223, 177, 239, 24, 247, 141, 251]
\end{lstlisting}

\paragraph{Height $m = 9$, pumping distance $d_{9} = 8$.}
\begin{lstlisting}[style=grid]
r=0  n=24  |S_1|=136
  A = [479, 369, 399, 184, 423, 365]
  B = [347, 338, 347, 365, 437, 149, 437, 365]
  C = [347, 338, 347, 365, 438, 219, 365, 437, 149, 247]
r=1  n=25  |S_1|=141
  A = [503, 285, 483, 58, 459, 365]
  B = [437, 149, 437, 365, 347, 338, 347, 365]
  C = [437, 149, 437, 365, 203, 442, 355, 349, 373, 269, 507]
r=2  n=26  |S_1|=148
  A = [479, 369, 399, 184, 423, 365]
  B = [347, 338, 347, 365, 437, 149, 437, 365]
  C = [347, 338, 347, 365, 438, 219, 365, 439, 152, 399, 369, 479]
r=3  n=27  |S_1|=153
  A = [479, 369, 399, 184, 423, 365]
  B = [347, 338, 347, 365, 437, 149, 437, 365]
  C = [347, 338, 347, 365, 437, 214, 347, 365, 439, 152, 399, 369, 479]
r=4  n=20  |S_1|=114
  A = [503, 285, 483, 58, 459, 365]
  B = [437, 149, 437, 365, 347, 338, 347, 365]
  C = [438, 219, 365, 437, 149, 247]
r=5  n=21  |S_1|=119
  A = [479, 369, 399, 184, 423, 365]
  B = [347, 338, 347, 365, 437, 149, 437, 365]
  C = [203, 442, 355, 349, 373, 269, 507]
r=6  n=22  |S_1|=126
  A = [503, 285, 483, 58, 459, 365]
  B = [437, 149, 437, 365, 347, 338, 347, 365]
  C = [438, 219, 365, 439, 152, 399, 369, 479]
r=7  n=23  |S_1|=131
  A = [503, 285, 483, 58, 459, 365]
  B = [437, 149, 437, 365, 347, 338, 347, 365]
  C = [437, 214, 347, 365, 439, 152, 399, 369, 479]
\end{lstlisting}

\paragraph{Height $m = 10$, pumping distance $d_{10} = 18$.}
\begin{lstlisting}[style=grid]
r=0  n=36  |S_1|=224
  A = [958, 739, 797, 309, 877, 731, 658, 731, 877, 439, 728, 879]
  B = [289, 879, 728, 951, 109, 987, 530, 987, 109, 949, 733, 835, 370, 795, 749, 695, 728, 879]
  C = [289, 879, 728, 695, 685, 955]
r=1  n=55  |S_1|=341
  A = [958, 739, 797, 309, 877, 731, 950, 109, 987]
  B = [530, 987, 109, 951, 728, 879, 289, 879, 728, 695, 749, 779, 314, 867, 733, 949, 109, 987]
  C = [530, 987, 109, 951, 728, 879, 289, 879, 728, 695, 749, 795, 310, 876, 728, 951, 109, 987, 530, 987, 109, 950, 731, 845, 373, 797, 737, 959]
r=2  n=38  |S_1|=236
  A = [958, 739, 797, 309, 877, 731, 694, 676, 694, 731, 877, 439, 728, 879]
  B = [289, 879, 728, 951, 109, 987, 530, 987, 109, 949, 733, 835, 370, 795, 749, 695, 728, 879]
  C = [289, 879, 728, 695, 685, 955]
r=3  n=57  |S_1|=353
  A = [958, 739, 797, 309, 877, 731, 950, 109, 987]
  B = [530, 987, 109, 951, 728, 879, 289, 879, 728, 695, 749, 779, 314, 867, 733, 949, 109, 987]
  C = [530, 987, 109, 951, 728, 879, 289, 879, 728, 695, 685, 699, 706, 891, 269, 887, 720, 951, 109, 987, 530, 987, 109, 950, 731, 845, 373, 797, 737, 959]
r=4  n=40  |S_1|=249
  A = [958, 739, 797, 309, 877, 731, 950, 109, 987]
  B = [530, 987, 109, 951, 728, 879, 289, 879, 728, 695, 749, 779, 314, 867, 733, 949, 109, 987]
  C = [530, 987, 109, 951, 728, 879, 289, 879, 728, 695, 749, 539, 1014]
r=5  n=77  |S_1|=475
  A = [958, 739, 797, 309, 877, 731, 950, 109, 987, 530, 987, 109, 951, 728, 879, 289, 879, 728, 695, 685, 699, 706, 891, 269, 885, 725, 949]
  B = [109, 987, 530, 987, 109, 951, 728, 879, 289, 879, 728, 695, 749, 779, 314, 867, 733, 949]
  C = [109, 987, 530, 987, 109, 951, 728, 879, 289, 879, 728, 695, 685, 699, 706, 891, 269, 887, 720, 951, 109, 987, 530, 987, 109, 950, 731, 845, 373, 797, 737, 959]
r=6  n=60  |S_1|=371
  A = [958, 739, 797, 309, 877, 731, 950, 109, 987, 530, 987, 109, 951, 728, 879, 289, 879, 728, 695, 685, 699, 706, 891, 269, 885, 725, 949]
  B = [109, 987, 530, 987, 109, 951, 728, 879, 289, 879, 728, 695, 749, 779, 314, 867, 733, 949]
  C = [109, 987, 530, 987, 109, 951, 728, 879, 289, 879, 728, 695, 749, 539, 1014]
r=7  n=43  |S_1|=268
  A = [503, 797, 739, 690, 731, 877, 439, 728, 879]
  B = [289, 879, 728, 951, 109, 987, 530, 987, 109, 949, 733, 835, 370, 795, 749, 695, 728, 879]
  C = [289, 879, 728, 951, 109, 987, 530, 987, 109, 950, 731, 845, 373, 797, 737, 959]
r=8  n=44  |S_1|=273
  A = [927, 689, 751, 792, 375, 845, 715, 442, 867]
  B = [733, 949, 109, 987, 530, 987, 109, 951, 728, 879, 289, 879, 728, 695, 749, 779, 314, 867]
  C = [733, 949, 109, 987, 530, 987, 109, 951, 728, 879, 289, 879, 728, 695, 749, 539, 1014]
r=9  n=63  |S_1|=390
  A = [503, 797, 739, 690, 731, 877, 439, 728, 879, 289, 879, 728, 951, 109, 987, 530, 987, 109, 949, 725, 885, 269, 891, 706, 699, 685, 695]
  B = [728, 879, 289, 879, 728, 951, 109, 987, 530, 987, 109, 949, 733, 835, 370, 795, 749, 695]
  C = [728, 879, 289, 879, 728, 951, 109, 987, 530, 987, 109, 950, 731, 845, 373, 797, 737, 959]
r=10  n=46  |S_1|=286
  A = [958, 739, 797, 309, 877, 731, 950, 109, 987]
  B = [530, 987, 109, 951, 728, 879, 289, 879, 728, 695, 749, 779, 314, 867, 733, 949, 109, 987]
  C = [530, 987, 109, 951, 728, 879, 289, 879, 728, 439, 877, 731, 658, 731, 877, 309, 797, 737, 959]
r=11  n=47  |S_1|=292
  A = [999, 565, 989, 99, 954, 715, 845, 374, 795]
  B = [749, 695, 728, 879, 289, 879, 728, 951, 109, 987, 530, 987, 109, 949, 733, 835, 370, 795]
  C = [749, 695, 728, 879, 289, 879, 728, 951, 109, 987, 530, 987, 109, 950, 731, 845, 373, 797, 737, 959]
r=12  n=66  |S_1|=408
  A = [958, 739, 797, 309, 877, 731, 950, 109, 987, 530, 987, 109, 951, 728, 879, 289, 879, 728, 695, 685, 699, 706, 891, 269, 885, 725, 949]
  B = [109, 987, 530, 987, 109, 951, 728, 879, 289, 879, 728, 695, 749, 779, 314, 867, 733, 949]
  C = [109, 987, 530, 987, 109, 951, 728, 879, 289, 879, 728, 439, 877, 731, 658, 731, 877, 309, 797, 737, 959]
r=13  n=49  |S_1|=304
  A = [503, 797, 739, 690, 731, 877, 293, 877, 731, 950, 109, 987]
  B = [530, 987, 109, 951, 728, 879, 289, 879, 728, 695, 749, 779, 314, 867, 733, 949, 109, 987]
  C = [530, 987, 109, 951, 728, 879, 289, 879, 728, 439, 877, 731, 658, 731, 877, 309, 797, 737, 959]
r=14  n=50  |S_1|=310
  A = [927, 689, 751, 792, 375, 845, 715, 442, 867]
  B = [733, 949, 109, 987, 530, 987, 109, 951, 728, 879, 289, 879, 728, 695, 749, 779, 314, 867]
  C = [733, 949, 109, 987, 530, 987, 109, 951, 728, 879, 289, 879, 728, 439, 877, 731, 658, 731, 877, 309, 797, 737, 959]
r=15  n=33  |S_1|=206
  A = [503, 797, 739, 690, 731, 877, 439, 728, 879]
  B = [289, 879, 728, 951, 109, 987, 530, 987, 109, 949, 733, 835, 370, 795, 749, 695, 728, 879]
  C = [289, 879, 728, 695, 685, 955]
r=16  n=52  |S_1|=322
  A = [503, 797, 739, 690, 731, 877, 439, 728, 879]
  B = [289, 879, 728, 951, 109, 987, 530, 987, 109, 949, 733, 835, 370, 795, 749, 695, 728, 879]
  C = [289, 879, 728, 951, 109, 987, 530, 987, 109, 949, 733, 867, 434, 219, 109, 951, 728, 879, 289, 879, 728, 695, 749, 539, 1014]
r=17  n=53  |S_1|=328
  A = [503, 797, 739, 690, 731, 877, 439, 728, 879, 289, 879, 728, 951, 109, 987, 530, 987, 109, 949, 725, 885, 269, 891, 706, 699, 685, 695]
  B = [728, 879, 289, 879, 728, 951, 109, 987, 530, 987, 109, 949, 733, 835, 370, 795, 749, 695]
  C = [728, 879, 289, 879, 728, 695, 685, 955]
\end{lstlisting}

\paragraph{Height $m = 11$, pumping distance $d_{11} = 8$.}
\begin{lstlisting}[style=grid]
r=0  n=56  |S_1|=382
  A = [1918, 1475, 1597, 613, 1757, 1459, 1902, 216, 951, 1645, 1499, 1298, 1499, 1645, 950, 219, 1901, 1463, 1752, 879, 433, 219]
  B = [1898, 1451, 1709, 695, 1752, 1135, 1969, 219]
  C = [1898, 1451, 1709, 695, 1752, 1135, 1969, 219, 1898, 1451, 1709, 695, 1752, 1132, 1974, 219, 1901, 1463, 1752, 879, 1457, 1759, 608, 1599, 1473, 1919]
r=1  n=57  |S_1|=388
  A = [1911, 1373, 1379, 1466, 1675, 749, 1590, 1499, 1389, 1463, 1240, 1903, 289, 1903, 1240, 1463, 1901, 219, 1974, 1132, 1752]
  B = [695, 1709, 1451, 1898, 219, 1969, 1135, 1752]
  C = [695, 1709, 1451, 1898, 219, 1969, 1135, 1752, 695, 1709, 1451, 1898, 219, 1969, 1135, 1752, 695, 1709, 1451, 1898, 219, 1969, 1135, 2008, 55, 2029, 1051, 2038]
r=2  n=58  |S_1|=395
  A = [1015, 1565, 1507, 1330, 1499, 1645, 951, 216, 1902, 1459, 1757, 581, 1757, 1459, 878, 1752, 1463, 1901, 219, 1974, 1132, 1752]
  B = [695, 1709, 1451, 1898, 219, 1969, 1135, 1752]
  C = [695, 1709, 1451, 1898, 219, 1969, 1135, 1752, 695, 1709, 1451, 1898, 219, 1969, 1135, 1752, 695, 1709, 1451, 1898, 219, 1969, 1135, 2008, 55, 2029, 1051, 2038]
r=3  n=59  |S_1|=401
  A = [1911, 1493, 1589, 749, 1675, 1466, 1379, 1501, 1589, 621, 1755, 1458, 1891, 221, 1973, 1133, 1755, 694, 1709, 1451, 1898, 219, 1969, 1135, 1752, 695, 1709, 1451, 1386, 1499, 1585, 623, 1752, 1463, 1901, 219, 1974, 1132, 1752]
  B = [695, 1709, 1451, 1898, 219, 1969, 1135, 1752]
  C = [695, 1709, 1451, 1898, 219, 1969, 1135, 2008, 55, 2029, 1051, 2038]
r=4  n=52  |S_1|=354
  A = [1911, 1373, 1379, 1466, 1675, 749, 1589, 1501, 1379, 1458, 1755, 877, 1463, 1240, 1903, 289, 1903, 1240, 1463, 1901, 219, 1974, 1132, 1752]
  B = [695, 1709, 1451, 1898, 219, 1969, 1135, 1752]
  C = [695, 1709, 1451, 1898, 219, 1969, 1135, 1752, 695, 1709, 1451, 1898, 219, 1969, 1135, 2008, 55, 2029, 1051, 2038]
r=5  n=53  |S_1|=361
  A = [1015, 1565, 1507, 1330, 1499, 1645, 951, 216, 1902, 1459, 1757, 581, 1757, 1459, 878, 1752, 1463, 1901, 219, 1974, 1132, 1752]
  B = [695, 1709, 1451, 1898, 219, 1969, 1135, 1752]
  C = [695, 1709, 1451, 1898, 219, 1969, 1135, 1752, 695, 1709, 1451, 1898, 219, 1969, 1135, 1240, 1463, 1901, 203, 442, 1891, 1245, 1975]
r=6  n=54  |S_1|=368
  A = [1918, 1475, 1597, 613, 1757, 1459, 1902, 216, 951, 1645, 1499, 1298, 1499, 1645, 950, 219, 1901, 1463, 1752, 879, 433, 219]
  B = [1898, 1451, 1709, 695, 1752, 1135, 1969, 219]
  C = [1898, 1451, 1709, 695, 1752, 1132, 1974, 219, 1901, 1463, 1752, 876, 438, 219, 1901, 1463, 1752, 879, 1457, 1759, 608, 1599, 1473, 1919]
r=7  n=55  |S_1|=375
  A = [1911, 1493, 1589, 749, 1675, 1466, 867, 1757, 1461, 1901, 219, 1974, 1060, 1974, 219, 1901, 1463, 1752, 879, 433, 219]
  B = [1898, 1451, 1709, 695, 1752, 1135, 1969, 219]
  C = [1898, 1451, 1709, 695, 1752, 1135, 1969, 219, 1898, 1451, 1709, 695, 1752, 1132, 1974, 219, 1901, 1463, 1752, 879, 1457, 1759, 608, 1599, 1473, 1919]
\end{lstlisting}

\end{document}